\documentclass{thesis}

\usepackage{float} 
\usepackage{hyperref} 
\usepackage{graphicx}
\usepackage{amssymb} 
\usepackage{CJKutf8} 
\newcommand{\jpni}{\begin{CJK}{UTF8}{min}ニ\end{CJK}}
\newcommand{\jpniu}{\textsuperscript{\jpni}}

\DeclareOldFontCommand{\sl}{\normalfont\slshape}{\@nomath\sl} 

\usepackage{enumitem}
\newlist{exercise}{enumerate}{5}
\setlist[exercise]{
	label*=\arabic*.,
	ref=\arabic*,
	before={
		\paragraph*{Let's pause and ponder}\ \\ (solutions in the appendix)
	} 
}
\let\ex\item

\newcommand{\thesistitle}{Introduction to Arnold's $J^+$-Invariant}
\newcommand{\thesisauthor}{Alexander Mai}
\newcommand{\submissiondate}{September 2022}

\title{\thesistitle}
\author{\thesisauthor}

\begin{document}

\pagestyle{thesistitlepage}

\begin{center}
	\vspace{2.5cm}
	\includegraphics[scale=0.37]{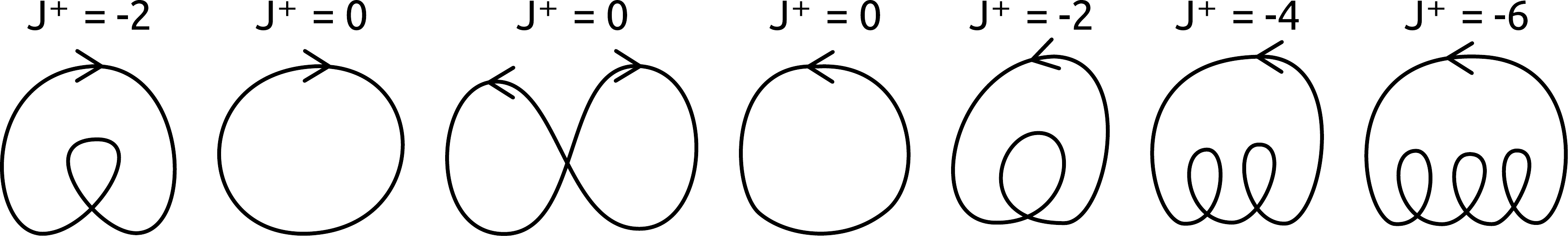}\\
	\vspace{0.9cm}
	\includegraphics[scale=0.34]{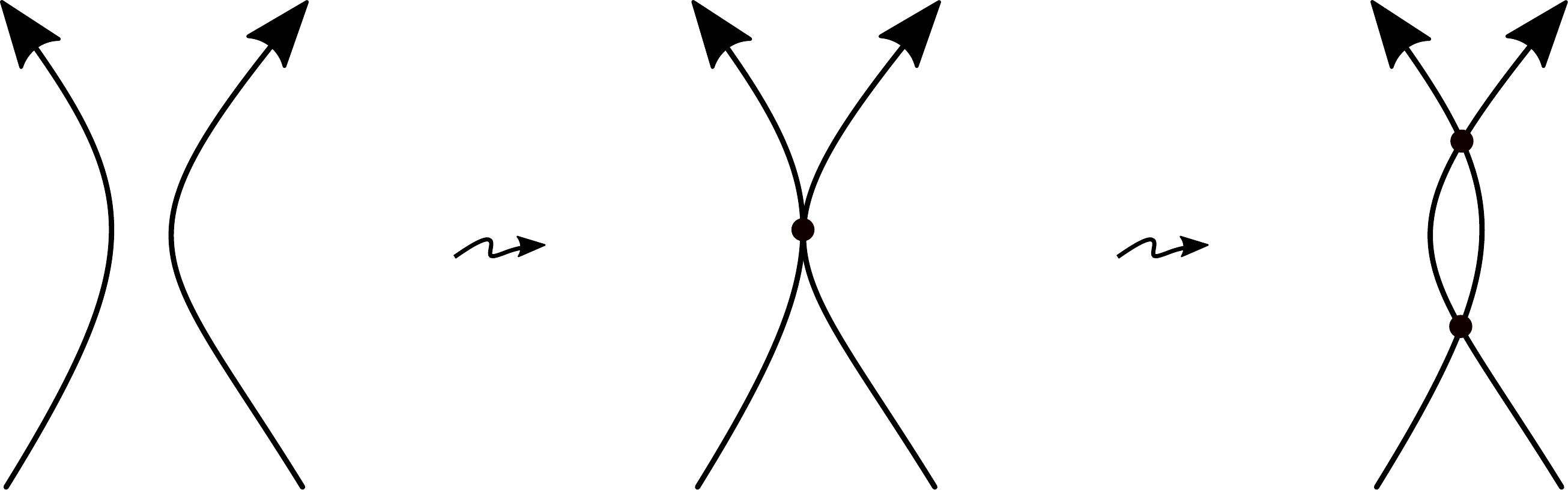}\\
	\vspace{0.5cm}
	\includegraphics[scale=0.27]{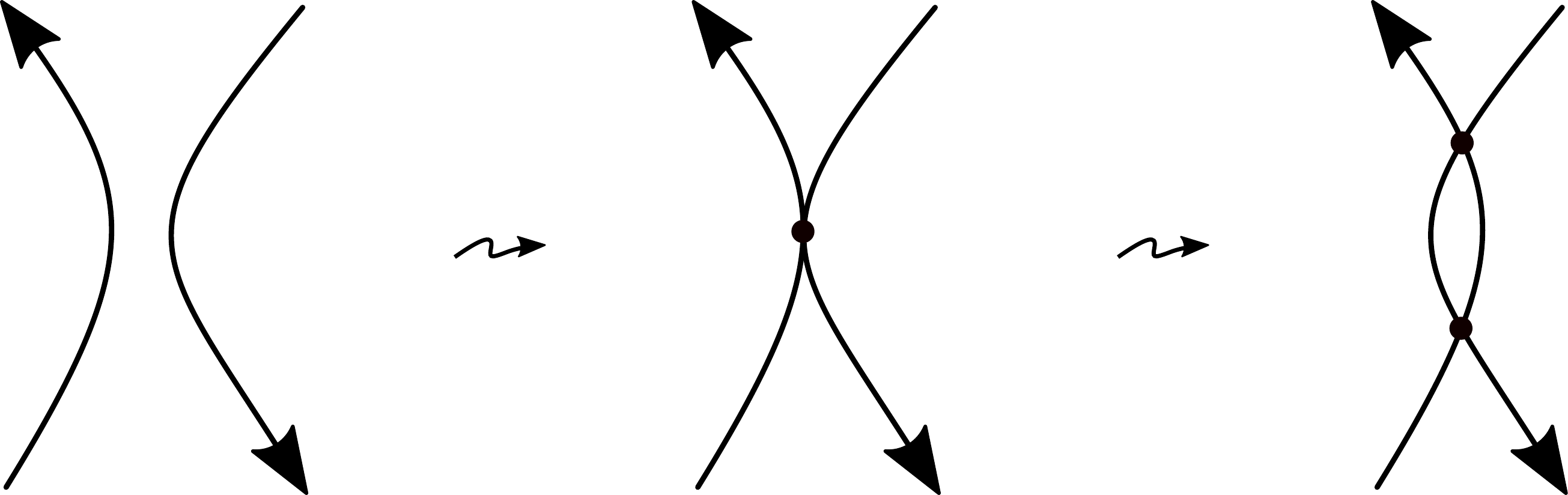}\\
	\vspace{0.5cm}
	\includegraphics[scale=0.24]{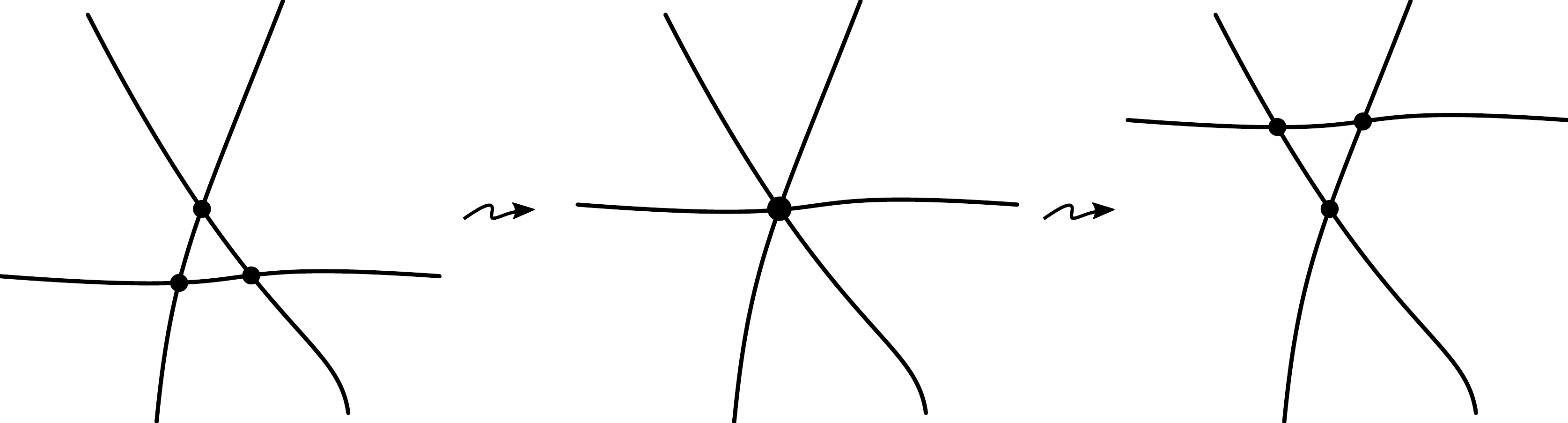}\\
	\vspace{2.5cm}
  	\normalsize
	{
		\centering
		\thesisauthor \\
		\href{mailto:alex.mai@posteo.net}{\Letter~alex.mai@posteo.net} \\
		\vspace{0.4cm}
		\submissiondate \\
		\vspace{0.4cm}
		University of Augsburg \\
	}
\end{center}

\newpage

\pagestyle{tocpage}
\tableofcontents

\newpage

\setcounter{page}{1}

\pagestyle{intro}

\section*{Abstract}
\addcontentsline{toc}{section}{\hspace{1.4em}Abstract}

\fancyhead[RO, LE]{Introduction}
\fancyhead[LO, RE]{}

We explore Arnold's $J^+$-invariant of immersions -- planar smooth closed curves with non-vanishing derivative, at most double points and only transverse intersections -- and computation methods like Viro's sum, among others.

Only basic undergraduate mathematics is needed to understand the contents of this introductory paper and everything we need that is above that is recalled or introduced. Examples, exercises and solutions are included for practice.

\section*{Introduction}
\addcontentsline{toc}{section}{\hspace{1.4em}Introduction}

In this introductory paper will thoroughly walk through the following: A regular immersed loop is a smooth map of a circle into the plane, i.e.\@ $q: S^1 \to \mathbb{C} ,$ with non-vanishing derivative. We identify the map with its image $K = q(S^1) \subset \mathbb{C} ,$ ignoring its parametrization and orientation. We call an immersed loop \emph{generic} if it has only transverse self-intersections and all of them are double points. For us all loops of interest are generic immersed loops, which we will simply call \emph{immersions} if not explicitly stated otherwise.

\emph{Vladimir Arnold} introduced three invariants for such immersions~\cite{arnold:paper}, of which the invariant $J^+$ is of special interest for this introductory paper. If we consider an immersion during a regular homotopy with only isolated non-generic moments, then the $J^+$-value of the immersion changes only under the first of the three so-called disasters -- direct self-tangencies, inverse self-tangencies, triple point intersections. If there are two immersions of the same rotation number, then by the \emph{Whitney--Graustein Theorem} one can be obtained from the other through a regular homotopy. Two such immersions have the same $J^+$-value if and only if during a regular homotopy from one to the other, the number of positive direct self-tangencies is the same as the number of negative direct self-tangencies. As a consequence, two immersions can only be homotopic to each other without direct self-tangencies -- a tangential double point where the directions agree -- if their $J^+$-value is the same.

Homotopies of immersions without direct self-tangencies are of interest for applications in planar orbital physics -- which are not discussed in this paper -- where the Hamiltonian is conserved as the sum of the potential and kinetic energy. For instance if the orbit of a particle moving in a conservative force field changes, direct self-tangencies cannot occur, but inverse self-tangencies can occur as long as velocity-dependent forces are present, like the Coriolis force or the Lorentz force. Without any velocity-dependent force neither direct nor inverse self-tangencies would be possible.

Other invariants based on $J^+$ can yield further applications, as for instance $\mathcal{J}_1$ and~$\mathcal{J}_2$ introduced by \emph{Cieliebak, Frauenfelder and van Koert} for any immersion $K \subset \mathbb{C}^* , \mathbb{C}^* \vcentcolon= \mathbb{C} \setminus \{ 0 \} ,$ see~\cite{kai:paper}. They are invariant under Stark--Zeeman homotopies (see~\cite{kai:paper}), which are of special interest for orbits of satellites in space but are not explored in this introductory paper.

After we finish the introduction to $J^+ ,$ as an extra, an approach to calculate~$J^+$ of interior sums (a notion to be introduced in this paper) of immersions is developed. We take any two immersions~$K, K' ,$ put $K'$ into a connected component $C$ of $\mathbb{C} \setminus K ,$ cross-connect the two immersions at two arcs where the orientations match, and call the resulting immersion~$K^{\times} .$ We then show that the following equation holds:
$$J^+(K^{\times}) = J^+(K) + J^+(K') -2 \cdot \omega_C(K) \cdot \operatorname{rot}(K') ,$$
with $\omega_C(K)$ the winding number of $K$ around $C$ and $\operatorname{rot}(K')$ the rotation number of $K' .$

In order to prove this, in Chapter~\ref{subsec:roteq} we find and prove an equation for the computation of the rotation number of an immersion~$K$ that only uses the winding numbers of the connected components of $\mathbb{C} \setminus K$ (denote them as $\Gamma_K$) and the index of each double point (denote them as $\mathcal{D}_K$), which we define to be the arithmetic mean of the winding numbers of the four connected components adjacent to a double point. We obtain the equation
$$ \operatorname{rot}(K) = \sum \limits_{C \in \Gamma_K} \omega_C(K) - \sum \limits_{p \in \mathcal{D}_K} \operatorname{ind}_p(K) .$$

First we recall some geometry basics in Chapter~\ref{sec:imbasics} and then introduce a few concepts of changing curves and events that can occur in Chapter~\ref{sec:homevents}.

Acknowledgements: This introductory paper would not exist without the successful completion of my bachelor thesis, which was only possible thanks to the extraordinarily patient guidance and motivation from Kai Cieliebak, Urs Frauenfelder, Ingo Blechschmidt, Julius Natrup, Florian Schilberth, Milan Zerbin, Leonie Nießeler and other friends from the University of Augsburg.

\newpage

\pagestyle{plain}

\section{Basics of immersions}
\label{sec:imbasics}

\begin{defi}[Immersed loop]
	\label{def:imloop}
	An \emph{immersed loop} is a regular loop, i.e.\@ a smooth map
	$$q: S^1 \to \mathbb{C}$$
	with non-vanishing derivative, which we will, by slight abuse of notation, from here on identify with its oriented image $K = q(S^1) \subset \mathbb{C} .$
\end{defi}

This abuse of notation is not dishonest, as we will not need any explicit definitions for immersed loops. Instead, we will draw many pictures -- in our heads and in this paper -- to talk about immersed loops. Note that orientation preserving reparametrizations of the immersion -- i.e.\@ changing how ``fast'' the curve is traversed without changing the direction -- have no effect on the oriented image and are consequently ignored from here on.

\begin{remark}
	In general, a \emph{path} in a topological space $X$ is a continuous map $f : [0, 1] \to X .$ It is called a closed path or a \emph{loop} if $f(0) = f(1) ,$ which is equivalent to defining a loop as a path with $S^1 ,$ the circle, as the domain of the path. A path (or loop) is called \emph{regular} if it is smooth and has non-vanishing derivative, so the derivative is non-zero everywhere. The image of a path (or loop) is often called a \emph{curve}. This is just a reminder -- all curves in this paper will be regular and closed, so we do not need these general definitions.
\end{remark}

Any curve that we can draw on paper, taking the following rules into consideration, is a valid visualization for an immersed loop:
\begin{itemize}
	\item ends where it started (loop),
	\item no interruptions (continuous),
	\item no edges (smooth, derivative never zero),
	\item is given one of two possible directions (oriented).
\end{itemize}

\begin{figure}[h!]
	\centering
	\includegraphics[scale=0.42]{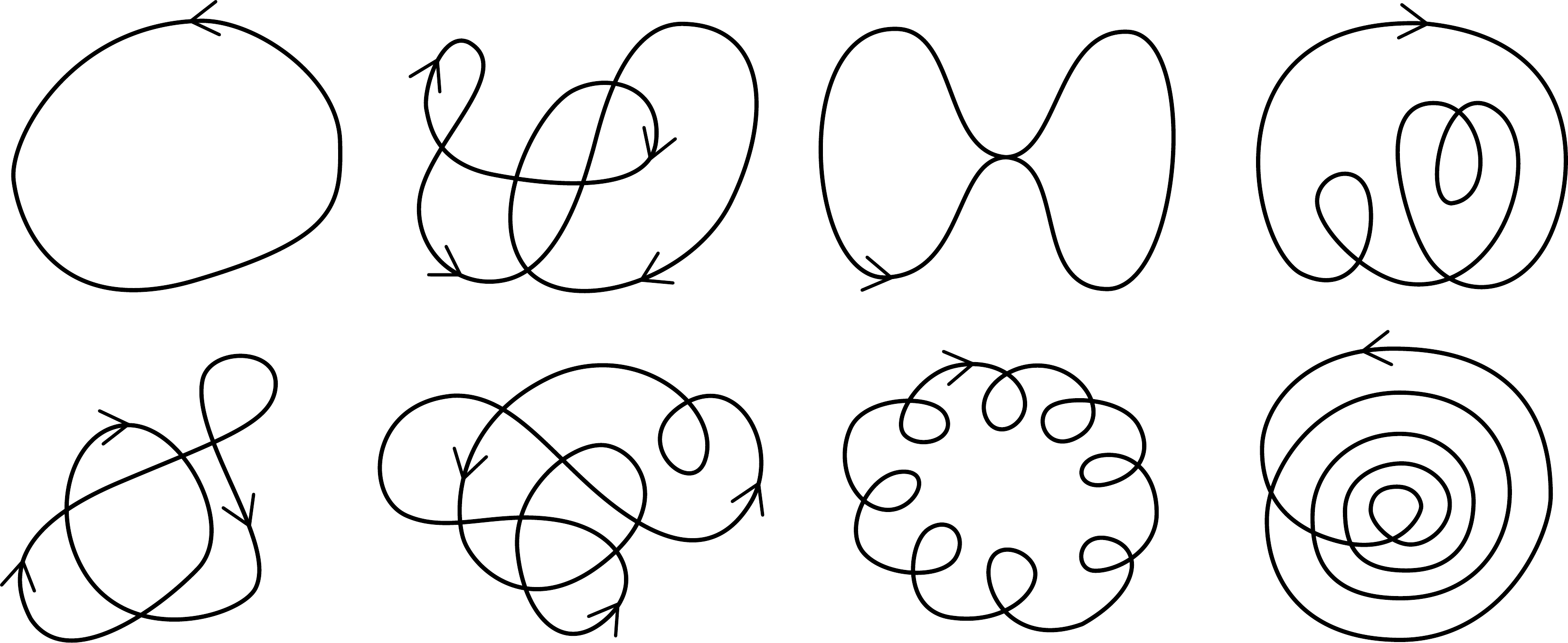}
	\caption{Examples of immersed loops.}
	\label{fig:imloopsex}
\end{figure}

The pictures in Figure~\ref{fig:imloopsex} represent immersed loops.

\begin{figure}[h!]
	\centering
	\includegraphics[scale=0.42]{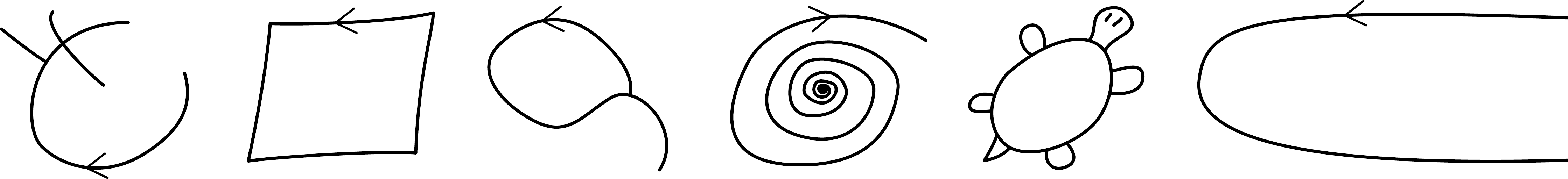}
	\caption{Not valid visualizations of immersed loops.}
	\label{fig:noimloopsex}
\end{figure}

The pictures in Figure~\ref{fig:noimloopsex} do not.

\begin{exercise}
    \ex Why is the third picture in the first row of Figure~\ref{fig:imloopsex} ambiguous? How many different immersions could be meant with the information of the picture and how could you fix it if you knew the immersion better?\label{ex:01}
    \ex Why are the pictures in Figure~\ref{fig:noimloopsex} not valid visualizations of immersed loops? Find at least one reason each.\label{ex:02}
	\ex Try to fix every picture in Figure~\ref{fig:noimloopsex} by adding some curves. Edges need some fantasy but are fairly easy.\label{ex:03}
\end{exercise}

If points are passed twice (thrice) by a curve, we call these points \emph{double points} (\emph{triple points}). Although a curve can pass any point an arbitrary amount of times, without loss of generality, we will not consider anything above triple points.

\begin{defi}[Generic immersion]
	\label{def:genim}
	We call an immersion a \emph{generic immersion} if it only has \emph{transverse} self-intersections which are at \emph{double points}.
\end{defi}

This is an important definition for our invariants, as $J^+$ (see Definition~\ref{def:jplus}) is only defined for generic immersed loops, i.e.\@ \emph{not} defined for immersed loops with tangential intersections or more-than-double points.

From here on we will often say \emph{immersions} to mean \emph{generic immersed loops} if not stated otherwise.

\begin{figure}[h!]
	\centering
	\includegraphics[scale=0.5]{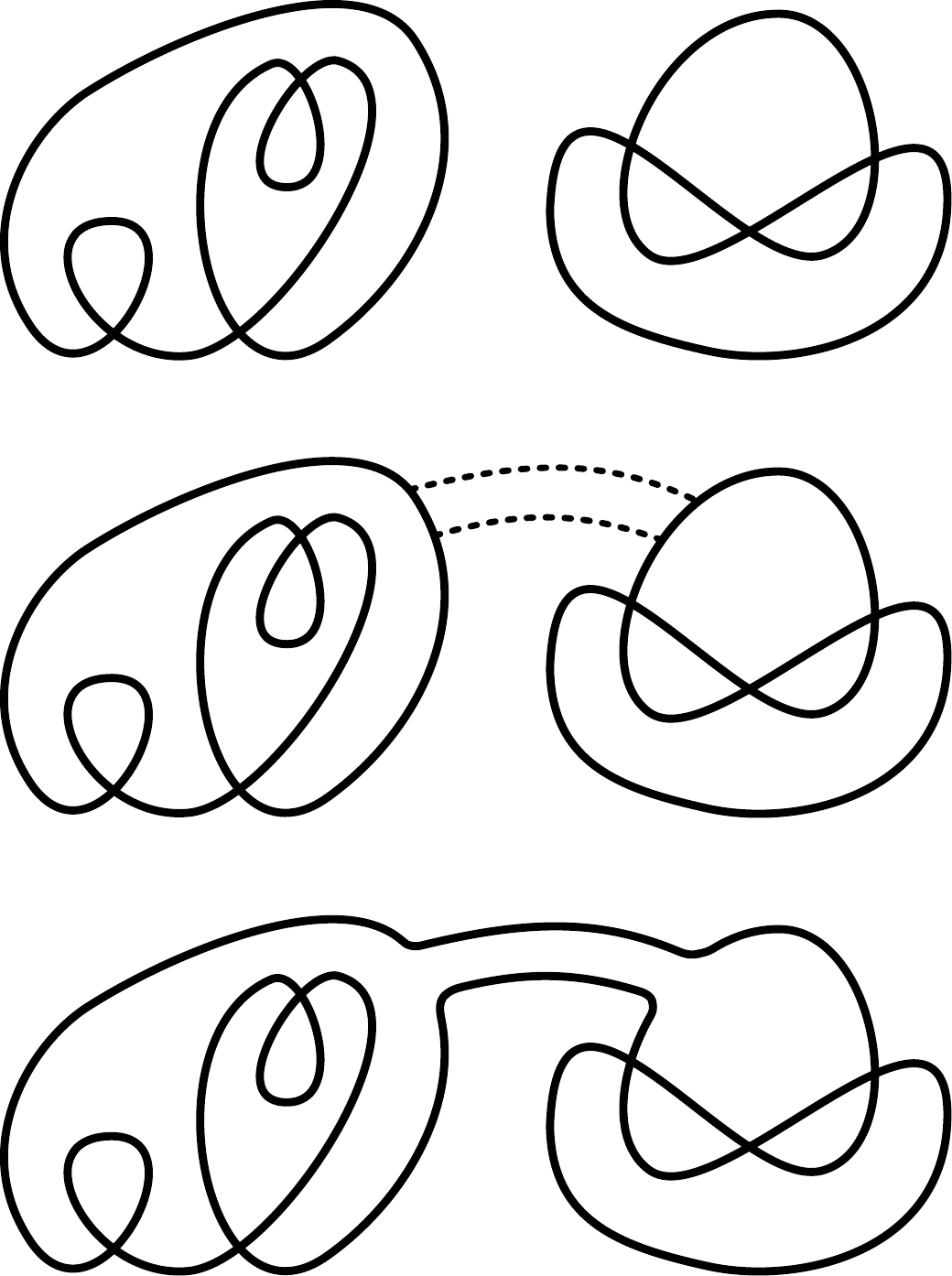}
	\caption{Connected sum of two immersions. We make sure to smoothen out the edges at the last step of the connection.}
	\label{fig:consum}
\end{figure}

\begin{defi}[Connected sum of immersions]
	\label{def:consum}
	The \emph{connected sum} of two immersions is created by joining them together as seen in Figure~\ref{fig:consum}.
\end{defi}

The immersions can be joined together at any outer arc. Of course this changes the resulting immersion.

Technically we should also be careful with the orientations of the immersions, which should not be changed by connecting them. There are two easy solutions to the problem of incompatible orientations at the connection:

\begin{itemize}
	\item put a sideways eight (later known as $K_0,$ see Figure~\ref{fig:standardcurves}) between the immersions, perform two connections, get a kind of \emph{cross-connected sum}
	\item ignore the problem and just change the orientation of one immersion.
\end{itemize}

Both ``solutions'' are fine for our purposes, as neither of them changes the invariant $J^+$ that we introduce below.

\begin{defi}[Winding number $\omega$]
	The \emph{winding number} $\omega_x(K)$ of an immersed loop $K$ is the amount of turns of the immersion around the point $x \in \mathbb{C} \setminus K,$ counting up (down) by one for every turn in counter-clockwise (clockwise) direction.
\end{defi}

For the sake of simplicity we forgo more formal definitions of the winding number.

Sometimes it is interesting to know the winding number of an immersion not around a single point, like the origin point, but around an arbitrary point within a \emph{connected component} of $\mathbb{C} \setminus K.$

\begin{defi}[Connected component]
	A \emph{connected component} $C$ is a maximal connected subset of a non-empty Euclidean space $X.$ An area $C$ is \emph{connected} if any two points can be connected by a curve within $C$ and \emph{maximal} if $\forall D \subseteq X, C \cap D \neq D {:} \, C \cup D$ is not connected.
\end{defi}

\begin{remark}
	For a generic immersed loop $K$ there are $n_K + 2$ many connected components, with $n_K$ the number of double points of $K.$ One of the components is unbounded, with winding number $0 .$
\end{remark}

We denote the winding number for all points in a connected component $C$ of $\mathbb{C} \setminus K$ as
$$\omega _C (K) \vcentcolon= \omega _y (K), y \in C \text{ arbitrary} .$$

\begin{figure}[h!]
	\centering
	\includegraphics[scale=0.4]{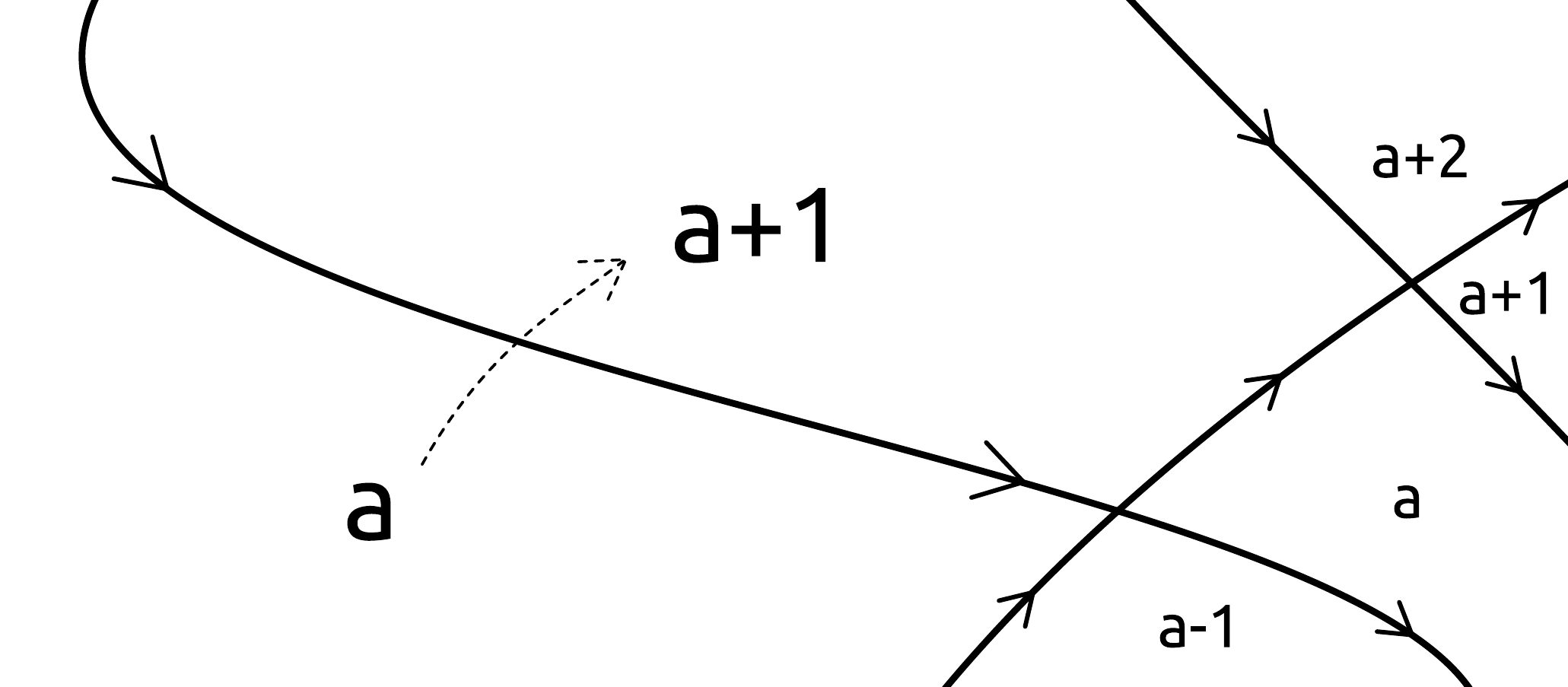}
	\caption{Changing winding number between adjacent connected components.}
	\label{fig:windingstep}
\end{figure}

The difference of the winding number of two connected components of $\mathbb{C} \setminus K$ that are adjacent to each other is always equal to $1.$ Looking at any connected component $C$ with a winding number of $a,$ the winding number of any adjacent component $C'$ is $a + 1$ ($a - 1$) if the immersion's arc, looking from $C$ to $C',$ is oriented from left to right (right to left). In other words, a flatlander living in one of the connected components can traverse an arc of the immersion to get to another connected component and the winding number of her location will increase (decrease) by $1$ if the traversed immersion arc is oriented from left to right (right to left), see Figure~\ref{fig:windingstep}.

With this trick it is easy to label all components' winding number in a picture of an immersion, starting from the unbounded component with winding number $0.$

\begin{thesislemma}
	\label{lem:dpwinding}
	The winding numbers around a double point $p$ are always of this form: let the lowest winding number around $p$ be $a \in \mathbb{Z}.$ Only one of the four components has winding number $a$, the component on the opposite site has winding number $a + 2$ and the other two have winding number~$a + 1 .$
\end{thesislemma}

\begin{figure}[h!]
	\centering
	\includegraphics[scale=0.4]{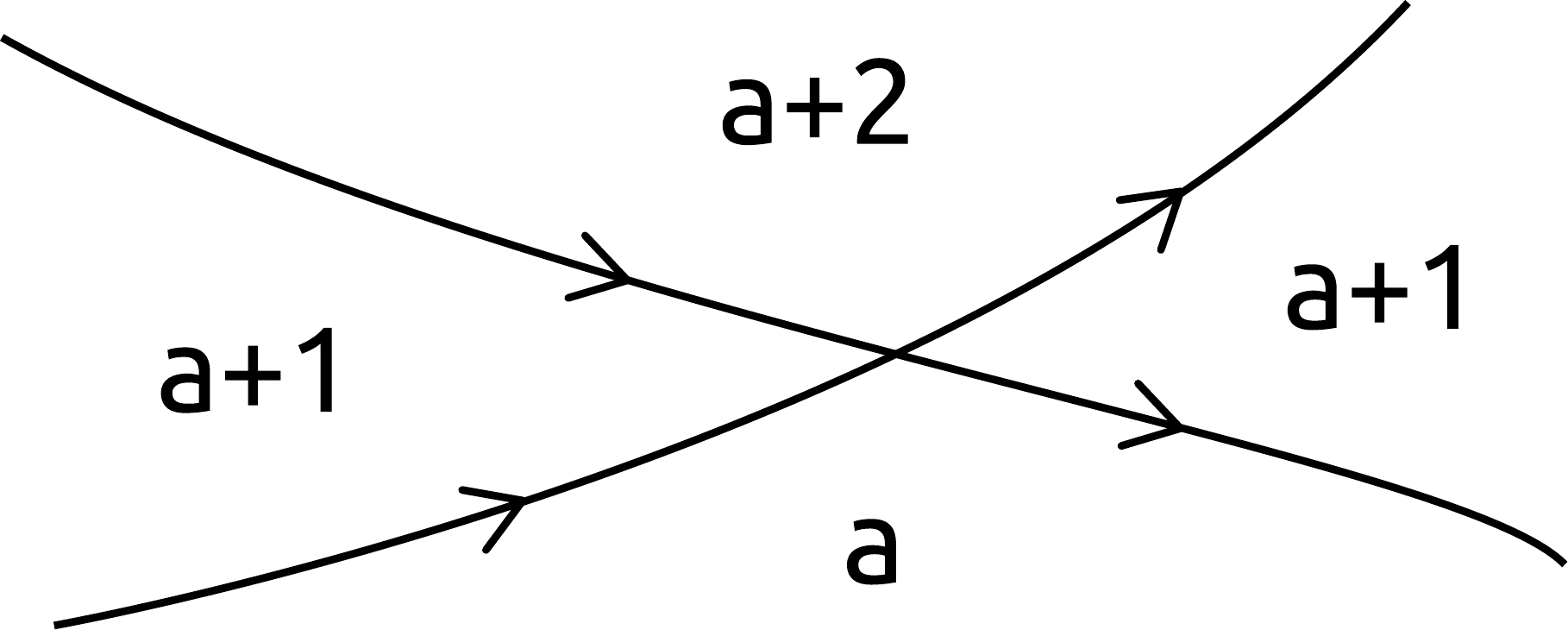}
	\caption{Winding numbers around a double point.}
	\label{fig:dpwinding}
\end{figure}

\begin{proof}
	See Figure~\ref{fig:dpwinding}. If $a$ is the lowest winding number around $p,$ the two adjacent components around $p$ have to be $a + 1,$ so the arcs to cross to get to these components have to be oriented left to right. This means the other two arcs also have to be oriented left to right and the last component has winding number $a + 2 .$
\end{proof}

\begin{defi}[Index $\operatorname{ind}_p(K)$]
	\label{def:indexdp}
	The \emph{index $\operatorname{ind}_p(K)$ of a double point $p$ in $K$} is the arithmetic mean of the winding numbers of the four connected components of $\mathbb{C} \setminus K$ adjacent to $p .$

	When calculating $\operatorname{ind}_p(K)$ and two corners of an adjacent component are equal to $p ,$ we count this one component and its winding number twice in the arithmetic mean to get a total of four adjacent components.
\end{defi}
\vspace*{0.4em}

\begin{remmm}
	\label{rem:indexdptwice}
	Using Lemma~\ref{lem:dpwinding} we can see that $\operatorname{ind}_p(K)$ is always equal to the winding number $a + 1$ appearing twice around $p,$ with $a \in \mathbb{Z}$ the lowest winding number around $p$ (same as in the lemma), as
	$$\dfrac{(a + (a + 1) + (a + 1) + (a + 2))}{4} = a + 1 .$$
\end{remmm}
\vspace*{0.4em}

\begin{remmm}
	\label{rem:indexdpside}
	It follows that whenever we draw a double point with the intersecting parts of the immersion both directed to the right (in the picture), then the index of that double point is equal to the winding number of the connected component to its right (or left). See Figure~\ref{fig:dpwinding}.

	This is also true if the intersecting parts are both directed from right to left.
\end{remmm}
\vspace*{0.4em}

\begin{defi}[Rotation number $\operatorname{rot}$]
	The \emph{rotation number} $\operatorname{rot}(q)$ of an immersed loop $q$ is equal to the winding number of $\dot{q},$ its derivative. So $\operatorname{rot}(q) = \omega_0(\dot{q}).$
\end{defi}

\begin{remark}
	Sometimes in other works the rotation number is called the \emph{(rotation) index} of an immersion.
\end{remark}

The rotation number is equal to how many times the tangent vector turns when travelling along the immersion once, counting up (down) by one for every turn in counter-clockwise (clockwise) direction.

\begin{remark}
	Changing the orientation of an immersion~$K$ changes the sign of its rotation number, the winding number of all the connected components in $\mathbb{C} \setminus K$ and the indices of all the double points of $K .$
\end{remark}

\begin{figure}[h!]
	\centering
	\includegraphics[scale=0.45]{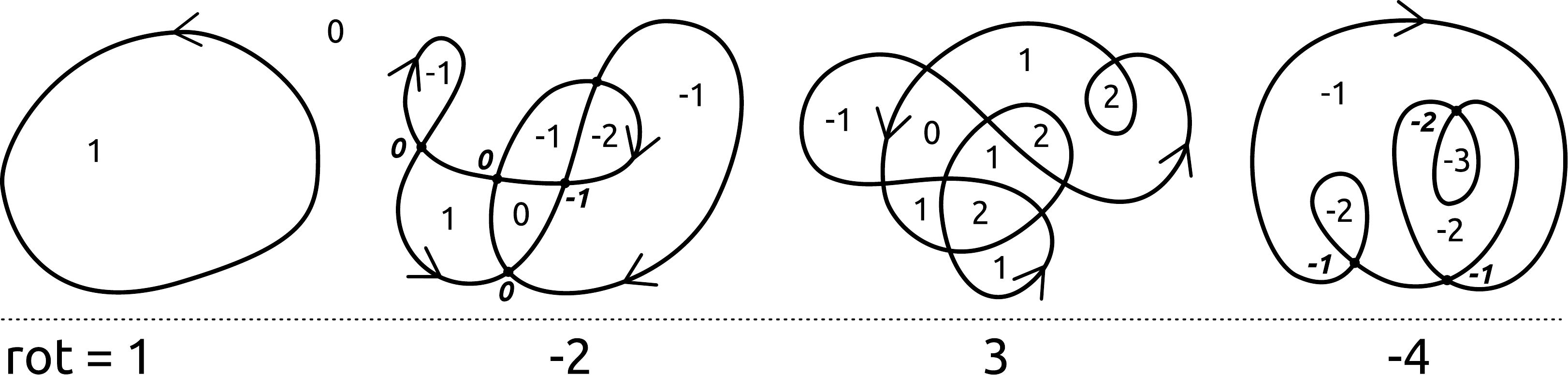}
	\caption{Immersed loops and their rotation number (below) and winding numbers for each connected component. The bold and cursive numbers are the indices of some double points.}
	\label{fig:exrotwindind}
\end{figure}

Figure~\ref{fig:exrotwindind} shows a few immersions and their winding numbers, double point indices and rotation number.

\begin{defi}[Interior loops]
	\label{def:intloops}
	An \emph{interior loop} within a connected component of an immersion~$K$ is the addition of a small loop at a boundary arc of the connected component. The loop is small enough to only intersect itself, not any other part of the immersion.

	A \emph{double interior loop} is an interior loop with an interior loop within itself. Analogous for \emph{triple}, \emph{quadruple} etc.\@ or simply \emph{$m$-interior loop}, $m \in \mathbb{N}^+ .$
\end{defi}

In Figure~\ref{fig:standardcurves} the first picture is a circle with one single interior loop. The last picture is a circle with three single interior loops. In Figure~\ref{fig:consum} the top left picture is a circle with a single interior loop on the left and a double interior loop on the right.

\newpage

\section{Events during homotopies of immersed loops}
\label{sec:homevents}

Looking at regular homotopies which \emph{do} pass through immersions with tangential intersections or triple points is fundamental to calculating the $J^+$-invariant of any non-trivial immersion.

\begin{defi}[Regular homotopy]
	A \emph{regular homotopy} between two immersed loops $q$ and $q'$ is a smooth map $$h : S^1 \times [0, 1] \to \mathbb{C},$$ with $h(\cdot, 0) = q$ and $h(\cdot, 1) = q'$ and $h(\cdot, t) : S^1 \to \mathbb{C}$ an immersed loop $\forall t \in [0, 1]$.

	We call two immersed loops \emph{regularly homotopic} if there is a regular homotopy between these two immersed loops.
\end{defi}

Remember that by the \emph{Whitney--Graustein Theorem} any two immersions are regularly homotopic if and only if they have the same rotation number.

Intuitively, a regular homotopy grabs one of the two immersions and deforms it into the other immersion and, most importantly, keeps it smooth with non-vanishing derivative -- i.e.\@ no edges -- during the deformation. A non-regular homotopy would have no such constraint that keeps the immersion smooth during the deformation.

Let us interpret the second variable of a regular homotopy as time passing going from one immersion to another, so we can talk about ``moments'' and a ``before'' and ``after''.

\begin{remark}
	The immersions during a regular homotopy do not have to always be generic. There can be $t_0 \in (0, 1)$ with $h(\cdot, t_0)$ not generic, i.e.\@ with self-tangencies or at least triple points. Without loss of generality, these moments $t_0$ are isolated within our regular homotopies and even in the non-generic case we have at most triple points. This is honest, because if the immersion during a homotopy is not generic in all but isolated moments, small perturbations of the homotopy can make the immersion generic in all but isolated moments.
\end{remark}

So during the process of deformation we can take note of three important isolated scenarios that can happen an arbitrary amount of times to the changing immersion: Direct self-tangencies, inverse self-tangencies, and triple points. In his paper, Arnold calls them ``perestroikas'' (Russian for \emph{restructuring, rearrangement, reorganisation}), which is how they are often called in literature since then. Others sometimes call them ``disasters''.

\vspace*{0.2em}
\begin{defi}[Direct (inverse) positive (negative) self-tangency]
	\label{def:selftang}
	A \emph{self-tangency} is the event of our immersion crossing itself tangentially. We call this self-tangency \emph{direct} (\emph{inverse}) if both parts involved in the crossing of the immersion are (are not) oriented in the same direction.

	We call a direct self-tangency \emph{positive} (\emph{negative}) if the number of double points increases (decreases) by $2 .$
\end{defi}
\vspace*{0.8em}

\begin{figure}[h!]
	\centering
	\includegraphics[scale=0.22]{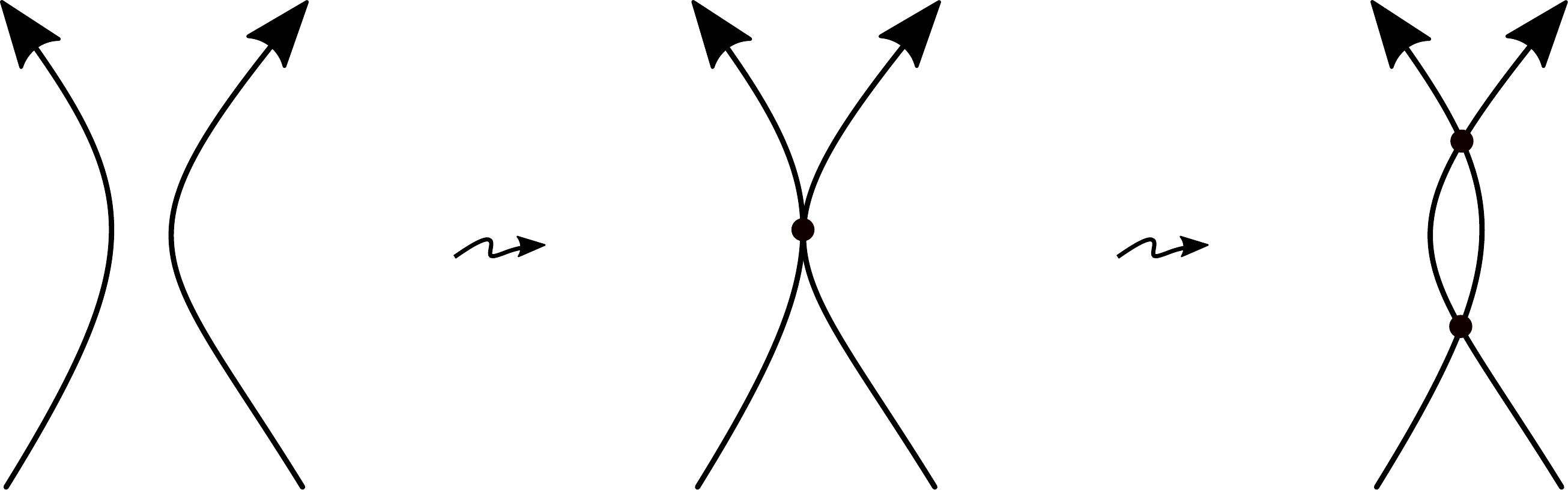}
	\caption{Direct self tangency during a regular homotopy. Left to right: positive. Right to left: negative.}
	\label{fig:jselfdirect}
\end{figure}

\begin{figure}[h!]
	\centering
	\includegraphics[scale=0.22]{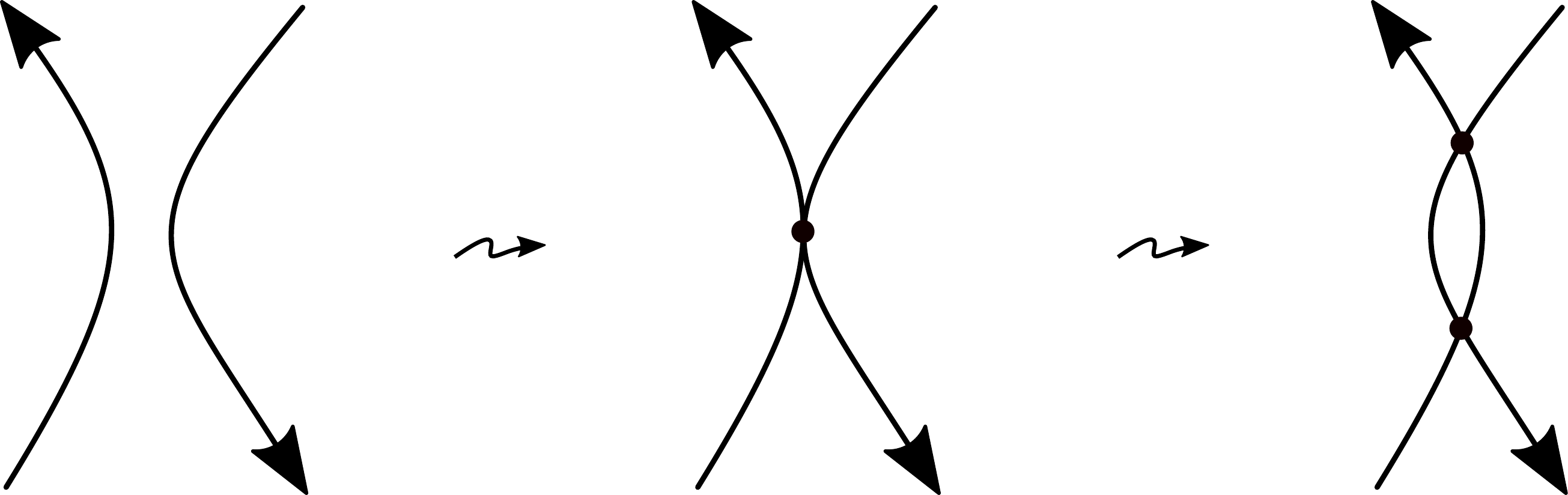}
	\caption{Inverse self tangency during a regular homotopy. Left to right: positive. Right to left: negative.}
	\label{fig:jselfinverse}
\end{figure}

In other words, if the self-tangency occurs at $t_0 \in (0, 1),$ then there is an $\varepsilon > 0$ so that the number of double points in $h(\cdot, t_0 - \varepsilon)$ and $h(\cdot, t_0 + \varepsilon)$ differs by two. See Figures~\ref{fig:jselfdirect} and~\ref{fig:jselfinverse}.

\begin{figure}[h!]
	\centering
	\includegraphics[scale=0.27]{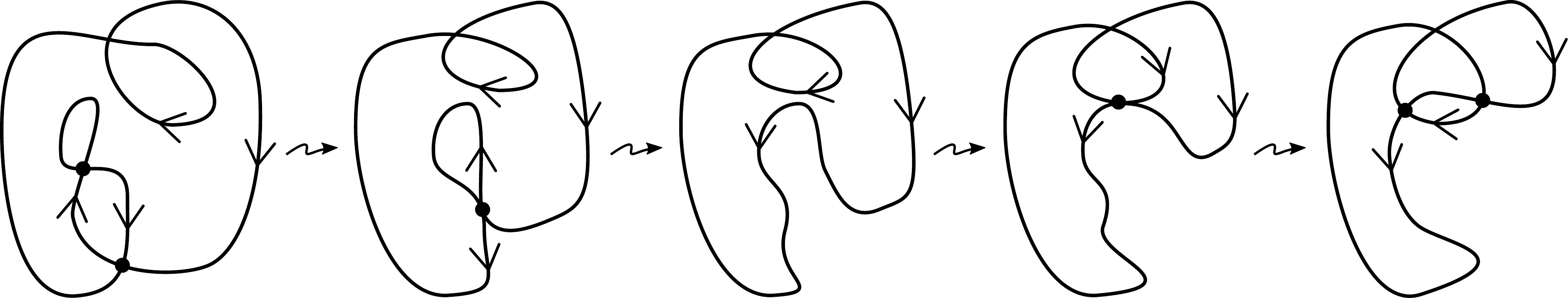}
	\caption{Regular homotopy between two immersed loops with (reading left to right) one negative inverse self-tangency (picture 2) and one positive direct self-tangency (picture 4), first removing, then adding two double points.}
	\label{fig:homtangex}
\end{figure}

Only direct self-tangencies will be important for our invariant $J^+$, inverse self-tangencies have no effect on it. We need to be careful not to confuse inverse/direct with positive/negative. One is about the orientation of the involved curve arcs. The other is about whether two double points appear or disappear.

\begin{figure}[h!]
	\centering
	\includegraphics[scale=0.22]{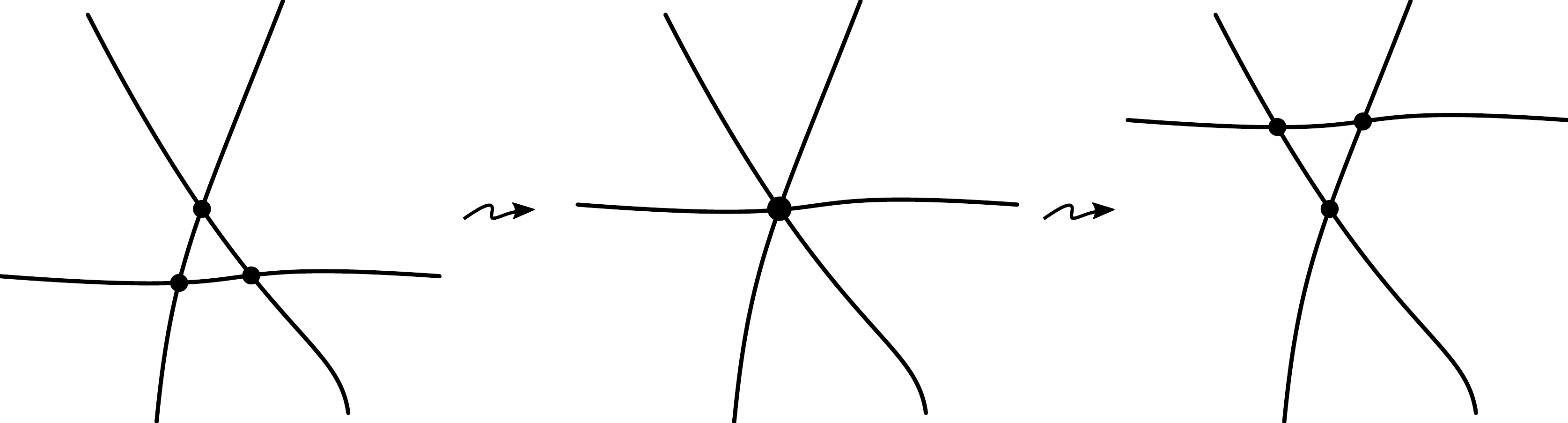}
	\caption{Triple point occurring during a regular homotopy.}
	\label{fig:jtrip}
\end{figure}

Another event that can occur, but has no effect on our invariant, is the crossing of triple points as seen in Figure~\ref{fig:jtrip}. The number of double points does not change at this event. Let us take note of that and be happy that we do not have to care about them for our intents and purposes.

\begin{figure}[h!]
	\centering
	\includegraphics[scale=0.3]{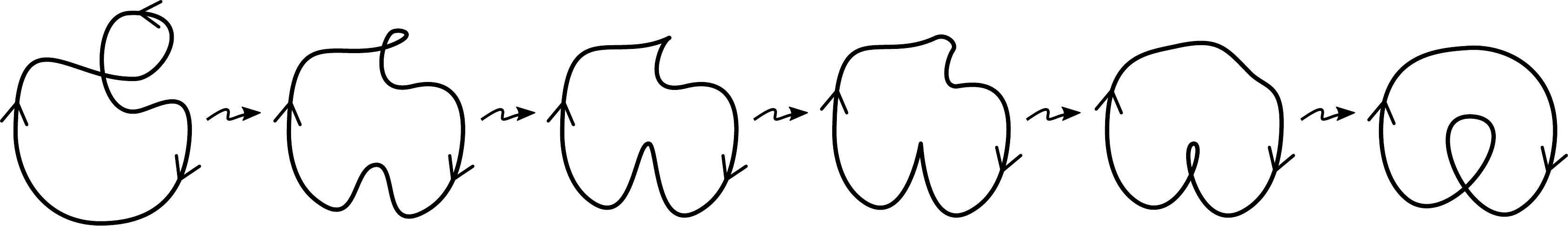}
	\caption{Example of cusps forming during a non-regular homotopy. Outer cusp in the third picture, inner cusp in the fourth picture.}
	\label{fig:imcusps}
\end{figure}

One last note is that in related literature, like in \cite{kai:paper}, other invariants based on $J^+$ are introduced that extend regular homotopies to allow for inner and outer cusps, which are then called \emph{Stark--Zeeman homotopies}. As cusps are not smooth, they will not appear in our homotopies, so be careful not to allow any when experimenting with homotopies, as they wipe out or create inner and outer loops, which would often change the invariant $J^+$. This also always changes the rotation number of the immersion. Figure~\ref{fig:imcusps} illustrates a homotopy with isolated cusps.

\newpage

\section{\texorpdfstring{The $J^+$-invariant}{The J+-invariant}}

In his paper \cite[Plane Curves, Their Invariants, Perestroikas and Classifications]{arnold:paper}, Arnold introduces the invariants $J^+, J^-$ and $St$ among others. His findings, especially on the $J^+$-invariant and its well-definedness, lay the foundations for the results of this paper.

\subsection{\texorpdfstring{Properties of $J^+$}{Properties of J+}}

\begin{figure}[h!]
	\centering
	\includegraphics[scale=0.4]{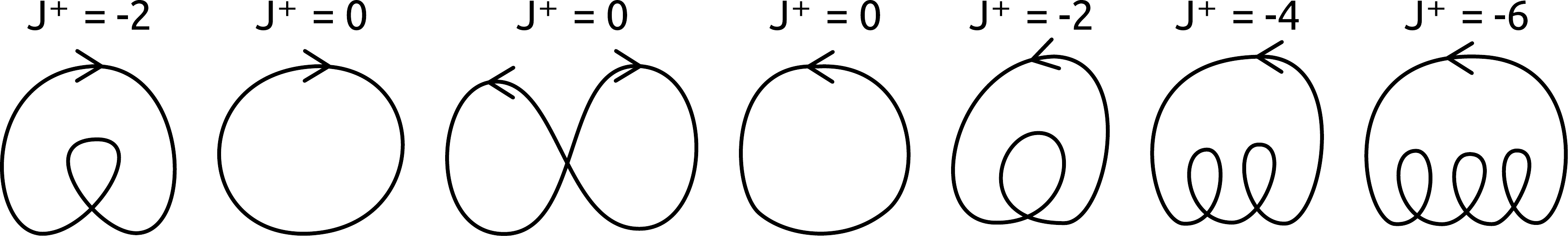}
	\caption{The standard curves. From left to right: $K_{-2}, K_{-1}, K_0, K_1, K_2, K_3, K_4 .$}
	\label{fig:standardcurves}
\end{figure}

Before diving into the technicalities of $J^+$, Figure~\ref{fig:standardcurves} are some fairly simple immersions and their respective values for $J^+$.

\begin{defi}[Standard curves $K_j$]
	\label{def:standardcurves}
	We call the immersions in Figure~\ref{fig:standardcurves} the \emph{standard curves} $K_j .$
	\begin{itemize}
		\item $K_0$ is the figure eight,
		\item $\forall j \neq 0: K_j$ is a circle with $|j| - 1$ many interior single loops that do not intersect with each other and with rotation number~$\operatorname{rot}(K_j) = j .$
	\end{itemize}
\end{defi}

\begin{defi}[$J^+$]
	\label{def:jplus}
	The invariant $J^+$ is a map
	$$\{ K \, | \, \text{$K$ is a generic immersed loop} \} \longrightarrow 2 \mathbb{Z},$$
	with
	$$K_0 \longmapsto 0$$
	and
	$$K_j \longmapsto -2 (|j| - 1) .$$

	Further, the invariant $J^+$
	\begin{itemize}\label{def:jplusprops}
		\item changes by $+2$ ($-2$) under positive (negative) direct self-tangencies, i.e.\@ tangent immersion crossings where the number of double points increase (decrease) and both involved tangent arcs have the same direction (see Definition~\ref{def:selftang}),
		\item is additive under \emph{connected sums} (see Figure~\ref{def:consum}),
		\item does not change under inverse self-tangencies or crossings through triple points,
		\item is independent of the orientation of the immersion.
	\end{itemize}
\end{defi}
\vspace*{0.8em}

In his paper, Arnold proves that with these properties $J^+$ is \emph{well-defined}. In his proof he shows that the difference of positive and negative direct self-tangency moments during a regular homotopy between two immersions is always the same, independent of the chosen regular homotopy \cite{arnold:paper}. This is important to remember when calculating $J^+$ of an immersion by showing that it is regularly homotopic to some other curve with known $J^+$, as any regular homotopy -- or the composition of many -- as simple or as complicated as it might be, will do.

Further, Arnold shows that there is in fact a unique invariant of generic immersions of fixed rotation number whose value does not change under inverse self-tangencies or triple point crossings, but increases (decreases) by a constant number $a_+$ under a positive (negative) direct self-tangency \cite[Theorem 2]{arnold:paper}. He then chooses $$a_+ = 2, \; \, K_0 \mapsto 0, \; \, K_j \mapsto -2(|j| - 1), \, \forall j \in \mathbb{Z} \setminus \{ 0 \}$$ and calls the unique invariant, with these choices, $J^+.$ The reasons for these choices are explained in his paper, one of which is allowing for additivity of connected sums.

Most of the time we only denote the orientation in the picture of an immersion to decide whether a self-tangency is direct or inverse, as $J^+$ is independent of the orientation of the immersion. With this in mind, some pictures of immersions from here on will lack any indication of orientation.

\subsection{\texorpdfstring{Calculation of $J^+$}{Calculation of J+}}

Let $K$ be an arbitrary generic immersed loop. By the \emph{Whitney--Graustein Theorem} we know that $K$ is regularly homotopic to the standard curve~$K_j$ (see Figure~\ref{fig:standardcurves}) with the same rotation number as $K,$ i.e.\@ to $K_{\operatorname{rot}(K)}.$

To calculate $J^+ (K),$ the most straight forward way is to find a regular homotopy between $K$ and the standard curve $K_{\operatorname{rot}(K)}$ and count how many positive direct self-tangencies $d_+$ and how many negative direct self-tangencies $d_-$ occur from $K$ to $K_{\operatorname{rot}(K)}$. By the properties of $J^+,$ see Definition~\ref{def:jplusprops}, it follows that:
$$J^+(K) = J^+(K_{\operatorname{rot}(K)}) - 2d_+ + 2d_-$$

Sometimes it is useful to use the additivity of connected sums to split an immersion up into two.

\begin{figure}[h!]
	\centering
	\includegraphics[scale=0.4]{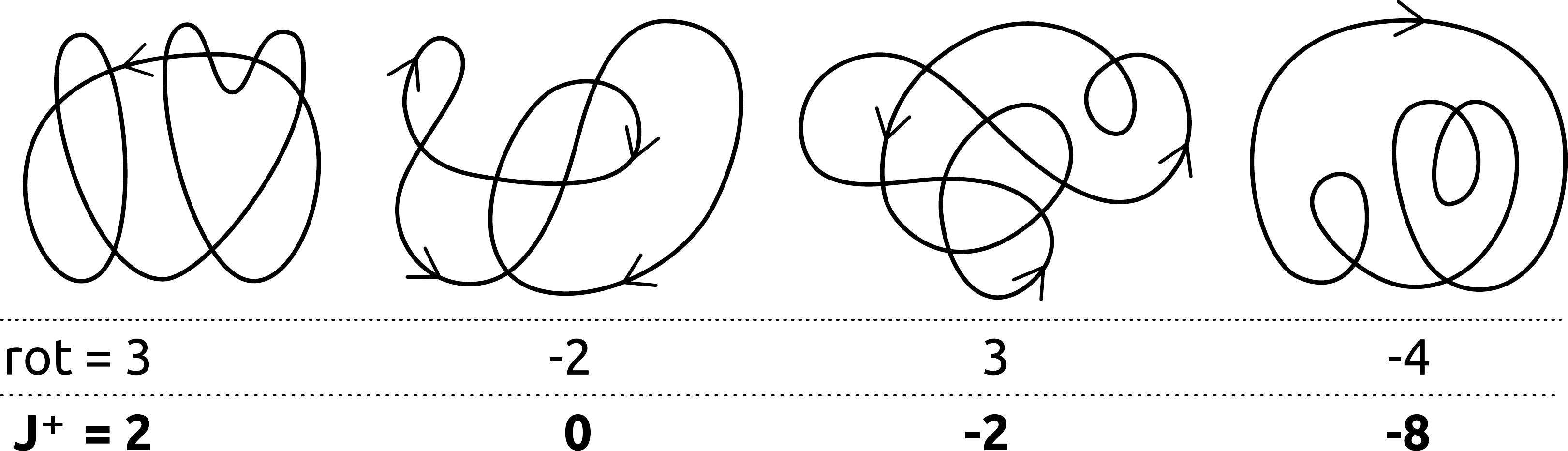}
	\caption{Some immersions and their $J^+$-value.}
	\label{fig:jexamples}
\end{figure}

\begin{exercise}[resume]
    \ex Try to verify the rotation number, winding numbers and double point indices for the immersions in Figure~\ref{fig:exrotwindind}.\label{ex:04}
	\ex Try to verify $J^+$ for the immersions in Figure~\ref{fig:jexamples}.\label{ex:05}
	\ex Calculate $J^+$ for the immersions at the top left, top right and bottom in Figure~\ref{fig:consum}.\label{ex:06}
	\ex How would $J^+$ of any $K_j, |j| > 2,$ change if neighboring interior loops would intersect each other a little, i.e.\@ creating two new double points? What if they intersect even more, creating two more new double points?\label{ex:07}
	\ex How would $J^+$ of any $K_j, |j| > 1,$ change if interior loops would intersect the top part of the circle?\label{ex:08}
	\ex The standard curve $K_2$ is a circle with a single interior loop. What is $J^+$ of a circle with a double interior loop? What about a triple interior loop?\label{ex:09}
\end{exercise}

Finding homotopies between an immersion and its corresponding standard curve can prove to be very arduous. Especially if we want to introduce even minor changes to a curve with already known~$J^+.$ Two important tools to calculate $J^+$ faster are given next.

\begin{itemize}
	\item A formula to directly calculate $J^+$ of an immersion with sums over its double points and connected components (Viro's formula).
	\item How two immersions with known $J^+$ can be combined by ``interior sums'' and what the $J^+$-value of the resulting immersion is.
\end{itemize}

\subsubsection{Viro's formula}%
\label{subsubsec:viro}%
Our first tool is \emph{Viro's Formula} for $J^+$ \cite[3.2.B Corollary]{viro:paper}. Given any generic immersed loop it calculates $J^+$ using the number of double points, the winding numbers of all components and the arithmetic mean of the winding numbers around each double point.

\begin{thesislemma}\label{lem:viro}[Viro's Formula]
	Let $K$ be a generic immersed loop. Then
	$$ J^+(K) = 1 + n_K - \sum \limits_{C \in \Gamma_K} (\omega_C(K))^2 + \sum \limits_{p \in \mathcal{D}_K} (\operatorname{ind}_p(K))^2 ,$$
	with
	\begin{itemize}
		\item $n_K$ the number of double points of $K$
		\item $\Gamma_K$ the connected components of $\mathbb{C} \setminus K$
		\item $\mathcal{D}_K$ the double points of $K$
		\item $\operatorname{ind}_p(K)$ the index of the double point $p$ in $K ,$ see Definition~\ref{def:indexdp}
	\end{itemize}
\end{thesislemma}

\begin{remark}
	In Viro's paper, the formula is stated as
	$$J^+(K) = 1 - \sum \limits_{C \in \Gamma_K} (\omega_C(K))^2 + \sum \limits_{p \in \mathcal{D}_K} (1 + (\operatorname{ind}_p(K))^2) ,$$
	which is equal to the way we wrote it.

	The formula can be proven the same way we proved Proposition~\ref{prop:rotequation}, which is different from the proof Oleg Viro did in his paper \cite{viro:paper}.
\end{remark}

\paragraph*{Example}

\begin{figure}[h!]
	\centering
	\includegraphics[scale=0.45]{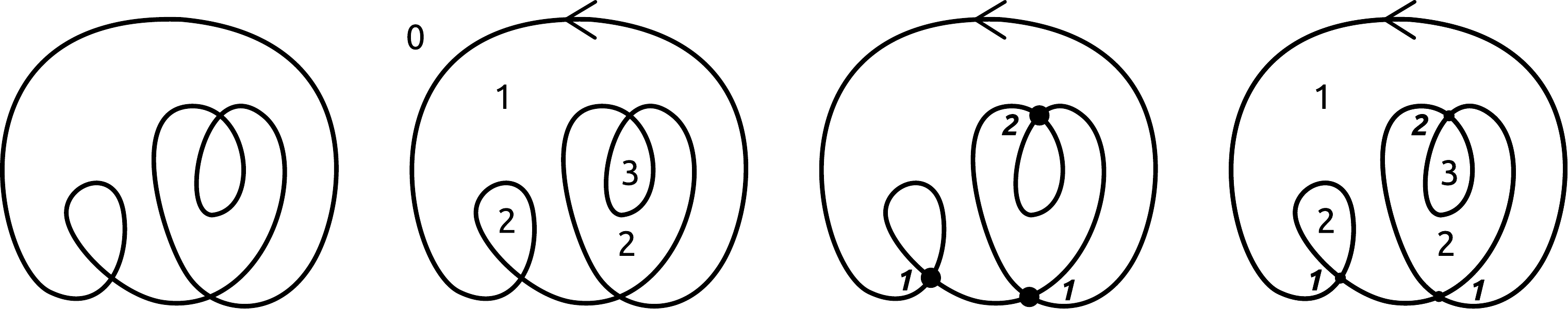}
	\caption{A circle with one single interior loop and one double interior loop. Picture one is just the immersion without orientation, picture two with winding numbers, picture three with double point indices, picture four with both.}
	\label{fig:viroex}
\end{figure}

Let us use Viro's formula to calculate $J^+$ of the immersion~$K$ from Figure~\ref{fig:viroex}. It is a circle with one single interior loop and one double interior loop.\label{example:viro}

In order to use Viro's formula, we need to count the number of double points, get the square of the winding number for all connected components of $\mathbb{C} \setminus K$ and the square of the double point indices.

The number of double points is $3 ,$ which can be easily seen in the third picture of Figure~\ref{fig:viroex}, where the double points are marked. So
$$ n_K = 3 . $$

Next we need all winding numbers. We started with an immersion~$K$ that has no orientation. For Viro's formula we will need the square of all winding numbers of the connected components, so it does not matter which orientation we choose. In this example we choose it so that the immersion has positive rotation number. With the orientation chosen, we can label all winding numbers of the connected components of~$\mathbb{C} \setminus K ,$ see the second picture of Figure~\ref{fig:viroex}. With this we can calculate
\begin{align*}
	\sum \limits_{C \in \Gamma_K} (\omega_C(K))^2
	&= 0^2 + 1^2 + 2^2 + 2^2 + 3^2 \\
	&= 0 + 1 + 4 + 4 + 9 \\
	&= 18 .
\end{align*}

Now the last thing we need is the double point indices. Again we only need the squares of them, so it is ok that we choose any of the two orientations. Remember that the index of a double point is the same as the winding number that appears twice around that double point (see Remark~\ref{rem:indexdptwice}). With this we can easily label all double point indices, see the third picture of Figure~\ref{fig:viroex}. With this we can calculate
\begin{align*}
	\sum \limits_{p \in \mathcal{D}_K} (\operatorname{ind}_p(K))^2
	&= 1^2 + 1^2 + 2^2 \\
	&= 1 + 1 + 4 \\
	&= 6 .
\end{align*}

We combine these observations to get $J^+(K)$ using Viro's formula:
\begin{align*}
	J^+(K)
	&= 1 + \underbrace{n_K}_{= 3} - \underbrace{\sum \limits_{C \in \Gamma_K} (\omega_C(K))^2}_{= 18} + \underbrace{\sum \limits_{p \in \mathcal{D}_K} (\operatorname{ind}_p(K))^2}_{= 6} \\
	&= 1 + 3 - 18 + 6 \\
	&= -8 .
\end{align*}

\subsection{Exercises}

In this subchapter we mainly take a look at a few different interesting immersions and how to calculate their $J^+$-value using homotopies and related tricks. It is not necessary to be proficient at calculating~$J^+$ of an immersion using homtopies in order to understand the results of this paper. This subchapter is mainly included because some experience in using homotopies to calculate $J^+$ of an immersion is an important skill that enables us to find new interesting results. Most of the results of this paper were envisioned drawing immersions and then applying homotopies in different ways. Once the idea for a result was there, it was simplified using Viro's formula.

Of course we can use Viro's formula to calculate $J^+$ all the time. Applying it is always the same procedure as it is just a matter of summing up integers, once the winding numbers are all labeled, see the example on page~\pageref{example:viro}. This makes it less prone to errors than using homotopies, as homotopies have to come to life in our heads and we need to find the isolated moments of direct self-tangencies during the homotopies. The main drawback is that using Viro's formula to see the change of $J^+$ after a homotopy is generally slower than checking for direct self-tangencies.

To avoid bloating this subchapter, we do not carry out the homotopies and calculations here, but instead display the solutions in the appendix at the end of this paper.

\begin{exercise}[resume]
    \ex Try to verify the rotation number and $J^+$ for the immersions in Figure~\ref{fig:interiortunnelsum}.\label{ex:10}

\begin{figure}[H]
	\centering
	\includegraphics[scale=0.38]{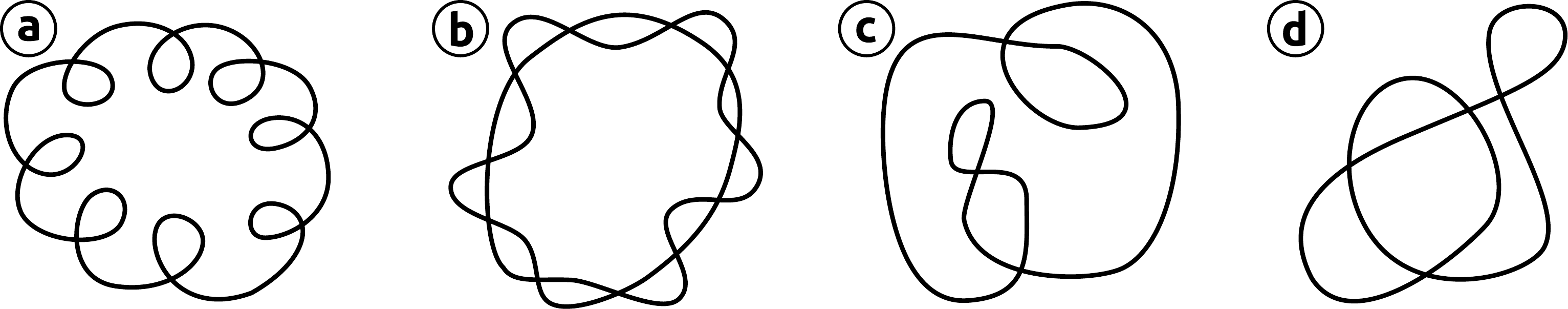}
	\caption{}
	\label{fig:manyjexamples1}
\end{figure}

\begin{figure}[H]
	\centering
	\includegraphics[scale=0.38]{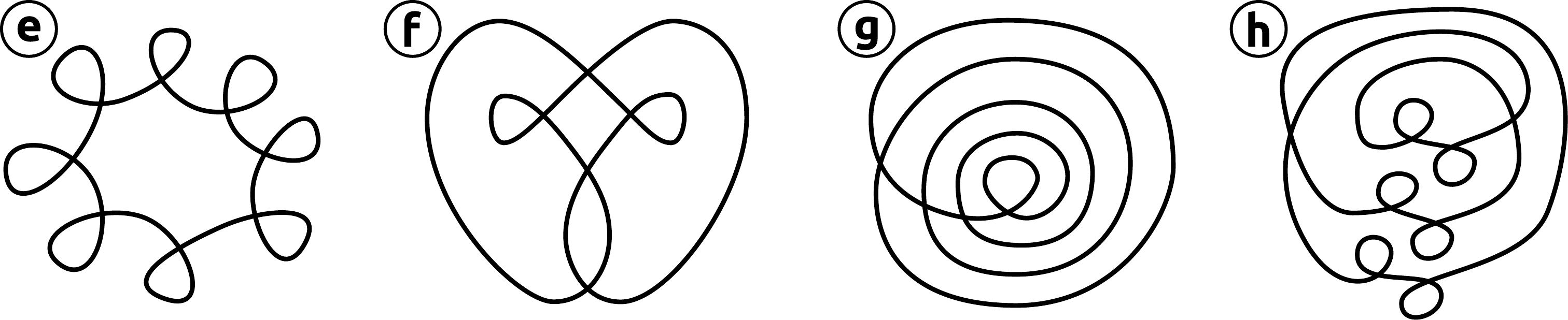}
	\caption{}
	\label{fig:manyjexamples2}
\end{figure}

\begin{figure}[H]
	\centering
	\includegraphics[scale=0.4]{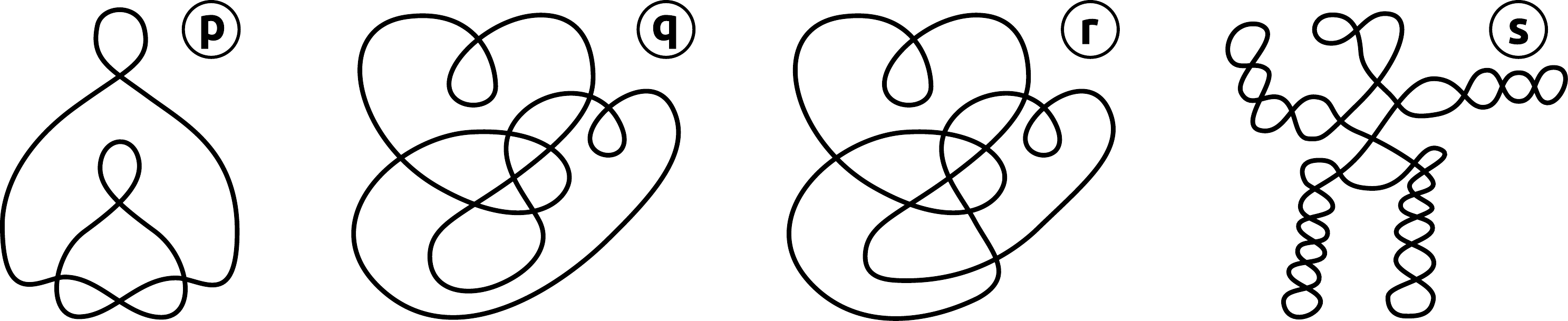}
	\caption{}
	\label{fig:manyjexamples3}
\end{figure}

\begin{figure}[H]
	\centering
	\includegraphics[scale=0.4]{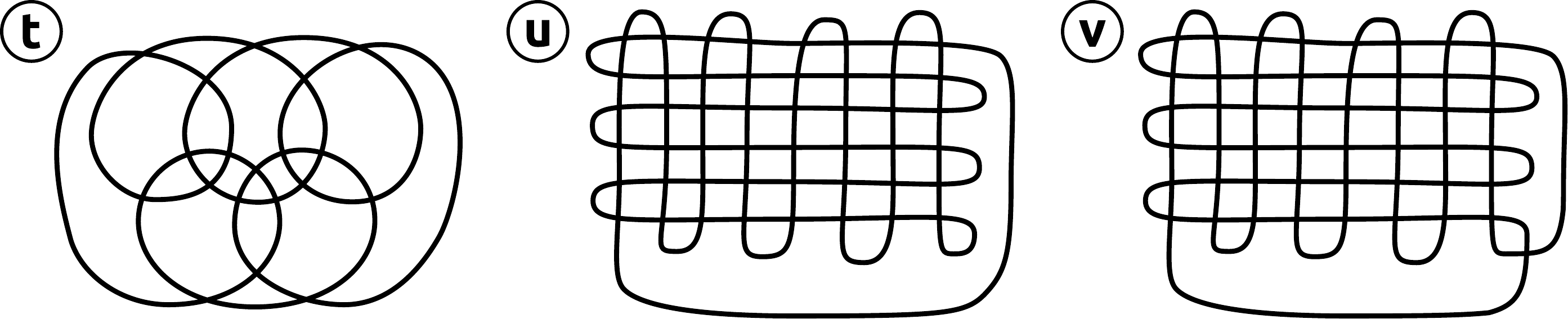}
	\caption{}
	\label{fig:manyjexamples4}
\end{figure}

	\ex Calculate $J^+$ for the immersions in Figures~\ref{fig:manyjexamples1}, \ref{fig:manyjexamples2}, \ref{fig:manyjexamples3} and~\ref{fig:manyjexamples4}.\label{ex:11}
	\ex Prove Viro's formula. Do not use Viro's formula or Theorem~\ref{th:interiorsum} [$J^+$ of interior sums].\label{ex:12}
\end{exercise}

\newpage

\section{Advanced Topic: Interior sum of immersions}
\label{sec:intsum}

This extra chapter is not part of the introduction to the $J^+$-invariant, but requires no additional knowledge and was added for readers who want to apply their knowledge on $J^+$ for practice or recreational reasons -- if not for using the contents of the following theorem and corollaries.

\begin{figure}[h!]
	\centering
	\includegraphics[scale=0.3]{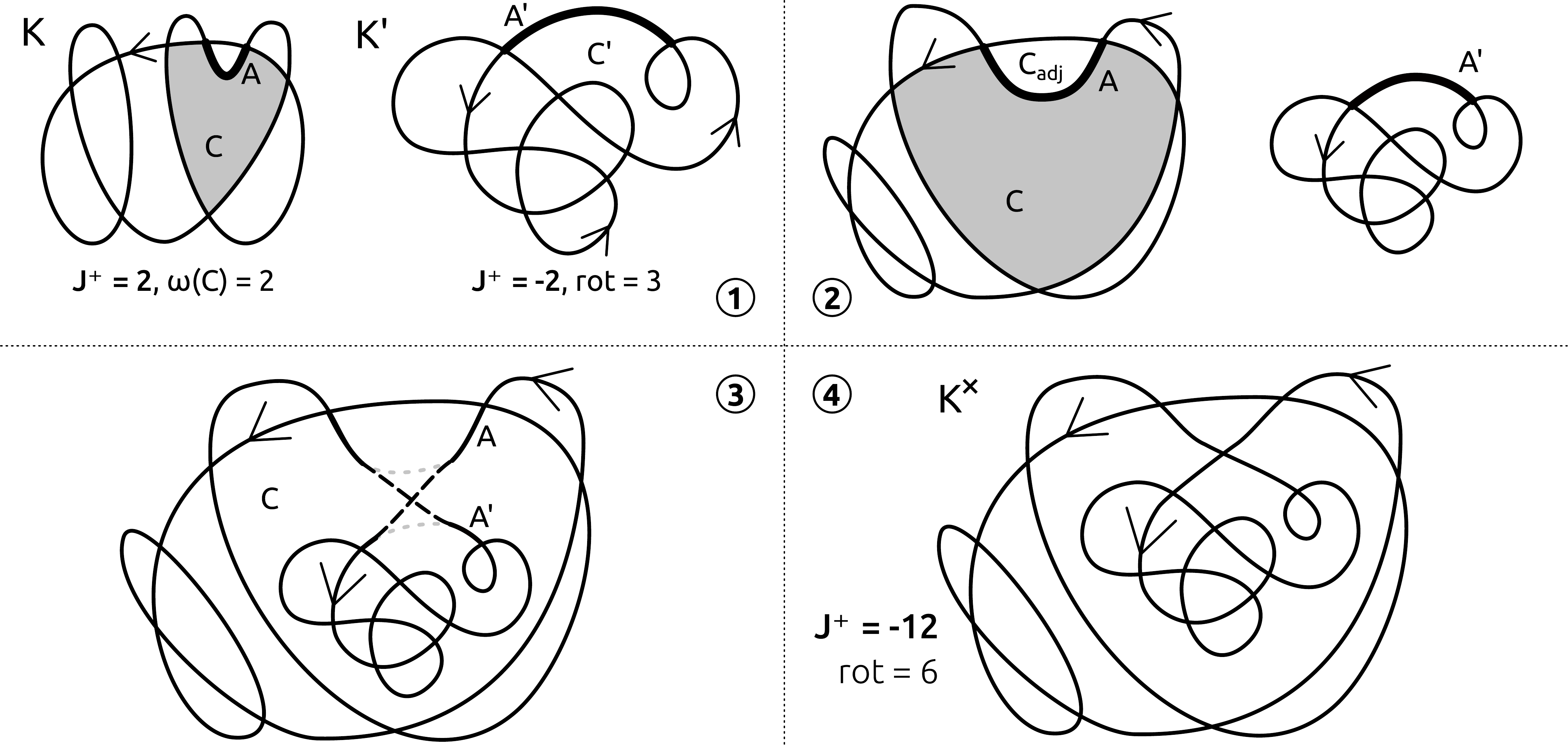}
	\caption{Interior sum of two immersions. The tinted connected component of immersion~$K$ is where we want to add $K'$ (first two pictures). The bold immersion arcs are where we want to connect the immersions together. In picture 2 we draw $K'$ smaller and the connected component $C$ bigger before putting them together in picture 3. There grey dotted lines hint at the former arcs $A$ and $A'$ which are now connected. Note that they are not connected with two parallel segments like in \emph{connected sums of immersions}, but instead with \emph{crossed} segments. The cross connection is visualized with dashed black lines. Picture 4 shows immersion $K^{\times}.$}
	\label{fig:interiorsum}
\end{figure}

\begin{defi}[Interior sum of immersions]
	\label{def:interiorsum}
	Let $K$ and $K'$ be two arbitrary immersions.

	Consider any bounded connected component $C$ of $\mathbb{C} \setminus K$ and any boundary arc~$A$ of~$C ,$ see Figure~\ref{fig:interiorsum}. Denote by $C_{\text{adj}}$ the unique connected component of $\mathbb{C} \setminus K$ that also has $A$ as a boundary arc and $C_{\text{adj}} \neq C .$

	Next consider any connected component $C'$ of $\mathbb{C} \setminus K'$ that shares a boundary arc $A'$ with the unbounded connected component.

	Denote the difference of winding numbers $\omega_C(K) - \omega_{C_{\text{adj}}}(K)$ as $\omega_{\text{adj}} .$

	If $\omega_{\text{adj}} \cdot \omega_{C'}(K') = -1 ,$ change the orientation of $K'$ so $\omega_{\text{adj}} = \omega_{C'}(K') \in \{ -1, 1\} .$ If the orientation of $K'$ needs to be changed, denote $K'$ with changed orientation still as $K'$ as if it was already equipped with the right orientation all along -- the interior sum is not defined for not matching orientations.

	Denote with $K^{\times}$ the immersion that we construct like this: put $K'$ into $C$ without intersections between $K$ and $K'$ -- it might need to be drawn smaller in visualizations -- and cross-connect the arcs $A$ and $A' .$

	We define this new immersion $K^{\times}$ as the \emph{interior sum of immersions $K, K'$ at boundary arcs $A, A'$ into the component $C$ of $\mathbb{C} \setminus K$}, with $A, A', C$ as above. See Figure~\ref{fig:interiorsum} for an illustration.
\end{defi}

The theorem presented and proven here is a byproduct of trying to prove one of the main results in a shortly to be published paper of mine about the change of $J^+$-like invariants at bifurcations, with inspiration from \cite[Lemma 4]{kai:paper}, which is stated in Corollary~\ref{cor:intloopkaiurs} and is a special case of the following theorem. This approach was later dropped and not used for anything at all. What we are left with is a fun theorem, a surprising formula for the rotation number (see the following Subchapter~\ref{subsec:roteq}) and a wealth of corollaries.

\begin{remark}
	\cite[Corollary 6.15]{kai:intsum} by \textnormal{Cieliebak, Frauenfelder and Zhao} can also be used to prove our Corollary~\ref{cor:interiortunnelsum}. It is proven differently than here and does not use Theorem~\ref{th:interiorsum}, but the theorem also directly follows from the corollary, so the results of this subchapter are not entirely new, but presented in a different light.
\end{remark}

\begin{theooo}[$J^+$ of interior sums]
	\label{th:interiorsum}
	If $K^{\times}$ is the interior sum of immersions $K, K'$ at boundary arcs $A, A'$ into the component $C$ of $\mathbb{C} \setminus K ,$ then
	\begin{equation*}
		J^+(K^{\times}) = J^+(K) + J^+(K') -2 \cdot \omega_C(K) \cdot \operatorname{rot}(K') .
	\end{equation*}
\end{theooo}

In order to prove this theorem we first need a surprisingly elegant formula to calculate the rotation number of an arbitrary immersion $K$ only from its connected components' winding number and double point indices.

\subsection{Rotation number from winding numbers}
\label{subsec:roteq}

\begin{propooo}[Rotation number from winding numbers]
	\label{prop:rotequation}
	Let $K$ be an arbitrary immersion -- i.e.\@ a generic regular closed curve. Then:
	\begin{equation}
		\label{eq:rotequation}
		\operatorname{rot}(K) = \sum \limits_{C \in \Gamma_K} \omega_C(K) - \sum \limits_{p \in \mathcal{D}_K} \operatorname{ind}_p(K)
	\end{equation}

	with
	\begin{itemize}
		\item $\Gamma_K$ the connected components of $\mathbb{C} \setminus K$
		\item $\mathcal{D}_K$ the double points of $K$
		\item $\operatorname{ind}_p(K)$ the index of the double point $p$ in $K ,$ see Definition~\ref{def:indexdp}
	\end{itemize}
\end{propooo}

We need the following lemma for one part of the proof of the proposition.

\begin{thesislemma}
	\label{lem:roteqsplit}
	A positive self-tangency of an immersion $K$ splits one connected component of $\mathbb{C} \setminus K$ into two and creates another new one, which increases the total amount of connected components by~two.
\end{thesislemma}
\vspace*{-0.75em}

\begin{proof}
	The completely new connected component is obvious and not the interesting part of the lemma. It is the one that is bounded only by the two arcs that were involved in the self-tangency (see Figure~\ref{fig:roteqinvt}, the one with winding number $a-1$).

	Now suppose that a positive self-tangency does not split a connected component~$C$ up into two separate connected components~$C_1$ and~$C_2 .$ So there is still a curve connecting all the points from $C_1$ with all points in $C_2 .$ This would mean that before the positive self-tangency there was a loop~$K^o$ in~$C$ that separates $\mathbb{C}$ into two disjoint areas and the immersion~$K$ had points in both areas without intersecting the loop~$K^o$, which is not possible for a continuous curve.
\end{proof}

\begin{proof}[Proof of Proposition~\ref{prop:rotequation}: Rotation number from winding numbers]
	We prove this in two easy steps. First we show that the equation is true for the standard curves $K_j$ (introduced later, see Definition~\ref{def:standardcurves}). Then we show that the right side of the equation does not change under regular homotopies.

	Let us prove the equation for the standard curves. We know that $\operatorname{rot}(K_j) = j$ by their definition.

	The standard curve $K_0$ has three connected components with winding numbers $0, -1, 1$ and one double point with index $0, $ see Figure~\ref{fig:roteqk0j}. So the equation is
	$$\operatorname{rot}(K_0) = (0 - 1 + 1) - 0 = 0 .$$

	The standard curve $K_j$ for $j > 0$ ($j < 0$) -- see Figure~\ref{fig:roteqk0j} -- has:
	\begin{itemize}
		\item the unbounded component with winding number $0$
		\item one connected component with winding number $1$ ($-1$)
		\item $|j| - 1$ many loops each bounding a connected component with winding number $2$ ($-2$)
		\item $|j| - 1$ many double points each with index $1$ ($-1$)
	\end{itemize}

	\begin{figure}[h!]
		\centering
		\includegraphics[scale=0.43]{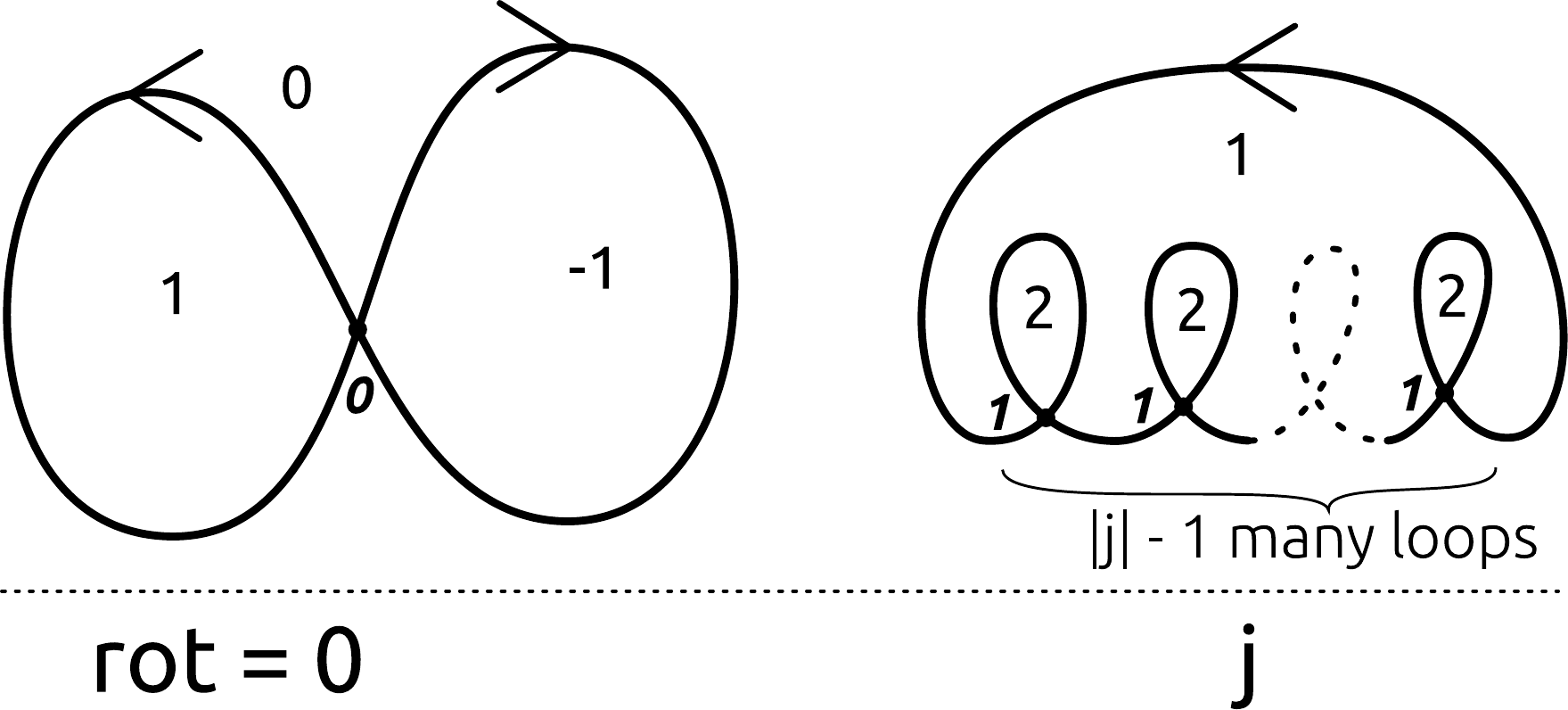}
		\caption{Standard curve $K_0$ on the left, $K_j$ on the right and their rotation number (below) and winding numbers for each connected component. The bold and cursive numbers are the indices of the double points. Here only for $j > 0$ multiply values by $-1$ for negative $j < 0 .$}
		\label{fig:roteqk0j}
	\end{figure}

	So the equation for $j > 0$ and for $j < 0$ is:
	$$ \operatorname{rot}(K_j) = \begin{cases}
		(0 + 1 + 2(j-1)) - 1(j-1) = j , &\text{if \( j > 0 \)} \\
		(0 - 1 - 2(-j-1)) + 1(-j-1) = j , &\text{if \( j < 0 \)}
	\end{cases} $$

	Next let us prove the invariance of $\sum \limits_{C \in \Gamma_K} \omega_C(K) - \sum \limits_{p \in \mathcal{D}_K} \operatorname{ind}_p(K)$ under regular homotopies for arbitrary immersions~$K .$

	Let $K$ be an arbitrary generic immersed loop. There are three events during regular homotopies that can change the topological properties of an immersion:
	\begin{itemize}
		\item direct self-tangency (see Figure~\ref{fig:jselfdirect})
		\item inverse self-tangency (see Figure~\ref{fig:jselfinverse})
		\item triple point crossing (see Figure~\ref{fig:jtrip})
	\end{itemize}

	For the first two events there is only one case to consider each. For the third event there are four. Each case is illustrated and verified to not change $\sum \limits_{C \in \Gamma_K} \omega_C(K) - \sum \limits_{p \in \mathcal{D}_K} \operatorname{ind}_p(K) .$

	In the pictures we denote some connected component to have winding number $a \in \mathbb{Z}$ and use the trick from Figure~\ref{fig:windingstep} to label the rest. Remember that the index of a double point is the same as the winding number that appears twice around that double point (see Remark~\ref{rem:indexdptwice}).

	When a direct self-tangency occurs, there are two new connected components (see Lemma~\ref{lem:roteqsplit}) with winding number $a+1$ and two new double points with index $a+1 ,$ see Figure~\ref{fig:roteqselft}, so $\sum \limits_{C \in \Gamma_K} \omega_C(K) - \sum \limits_{p \in \mathcal{D}_K} \operatorname{ind}_p(K)$ changes by
	$$2(a+1) - 2(a+1) = 0 .$$

	\begin{figure}[h!]
		\centering
		\includegraphics[scale=0.4]{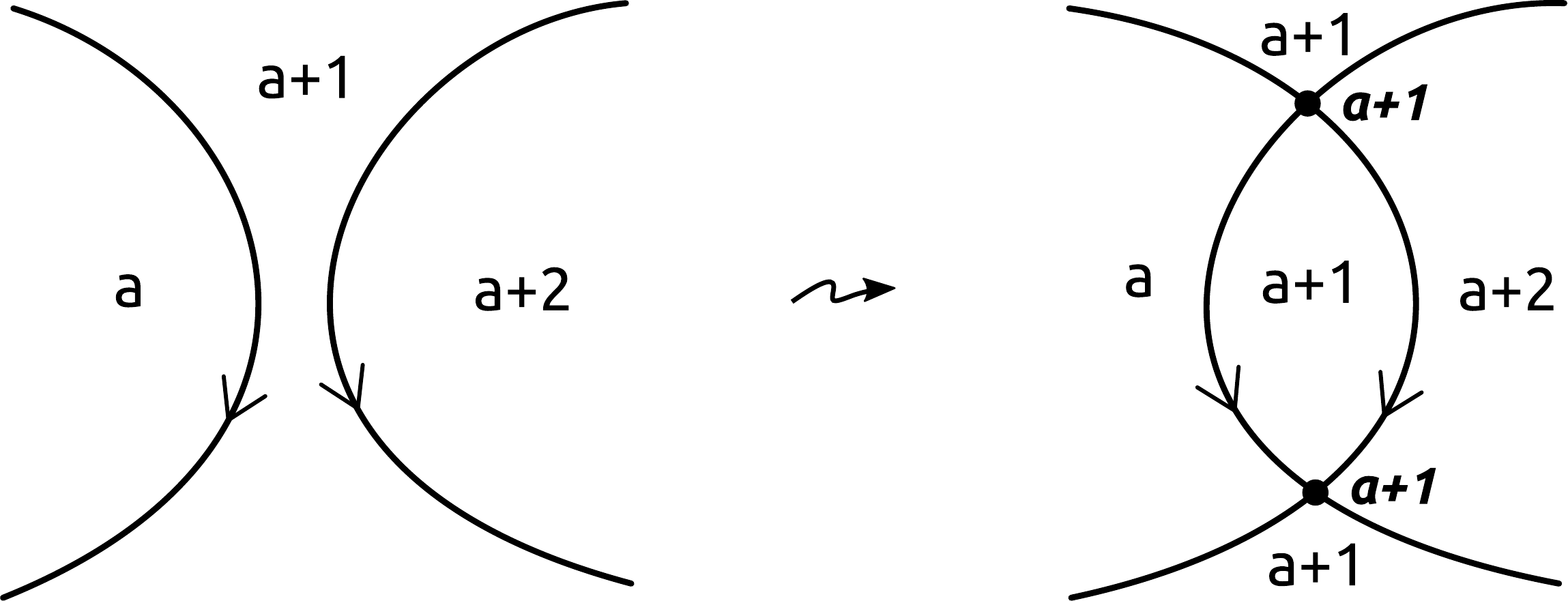}
		\caption{Winding numbers for each connected component and indices of the double points (bold and cursive) for a direct self-tangency.}
		\label{fig:roteqselft}
	\end{figure}

	When an inverse self-tangency occurs, there are two new connected components (see Lemma~\ref{lem:roteqsplit}) with winding number $a-1$ and $a+1$ and two new double points with index $a ,$ see Figure~\ref{fig:roteqinvt}, so $\sum \limits_{C \in \Gamma_K} \omega_C(K) - \sum \limits_{p \in \mathcal{D}_K} \operatorname{ind}_p(K)$ changes by
	$$((a-1) + (a+1)) - 2a = 0 .$$

	\begin{figure}[h!]
		\centering
		\includegraphics[scale=0.4]{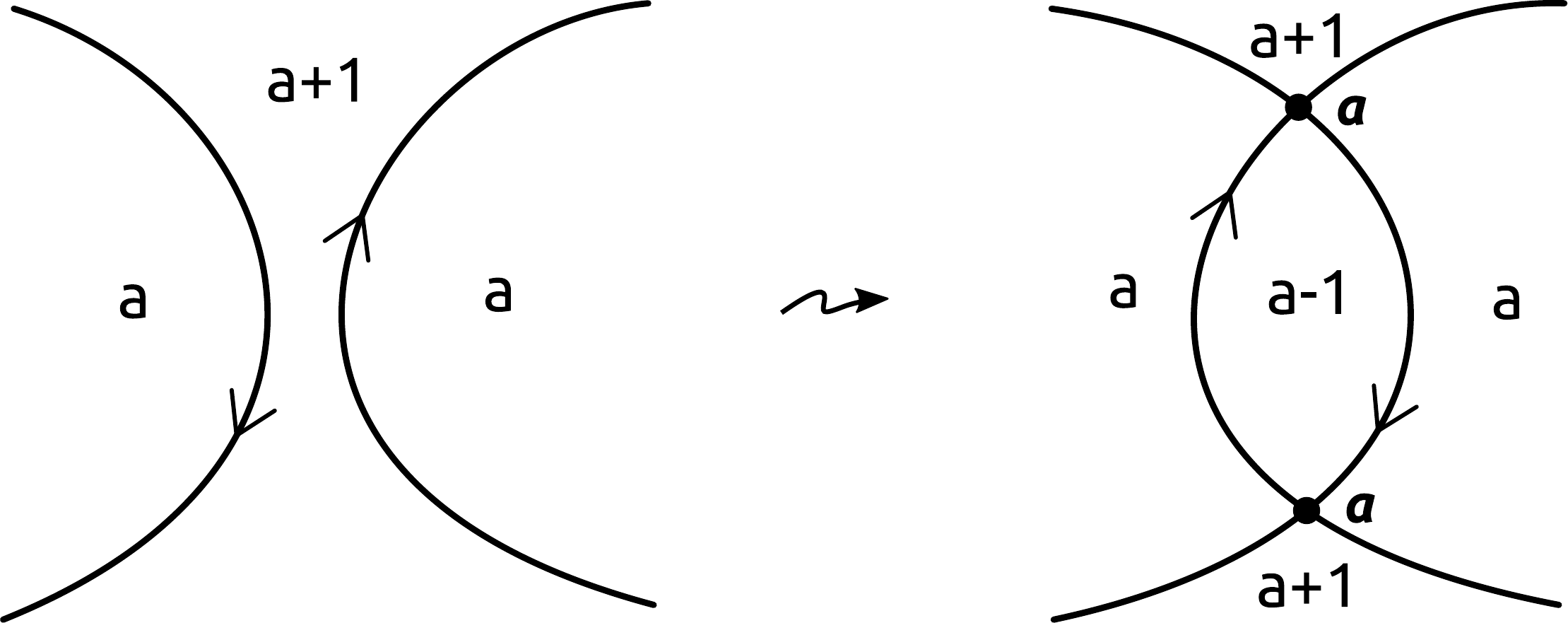}
		\caption{Winding numbers for each connected component and indices of the double points (bold and cursive) for an inverse self-tangency.}
		\label{fig:roteqinvt}
	\end{figure}

	When a triple point crossing occurs, one connected component is removed, a new one is added and three double points are changed. To keep the pictures easy to read, we pick some winding number to be $0$ instead of $a$ and go from there, but of course any $a \in \mathbb{Z}$ can be added to the pictures' winding numbers. This does not change the result, as all the $a$ cancel out.

	\begin{figure}[h!]
		\centering
		\includegraphics[scale=0.3]{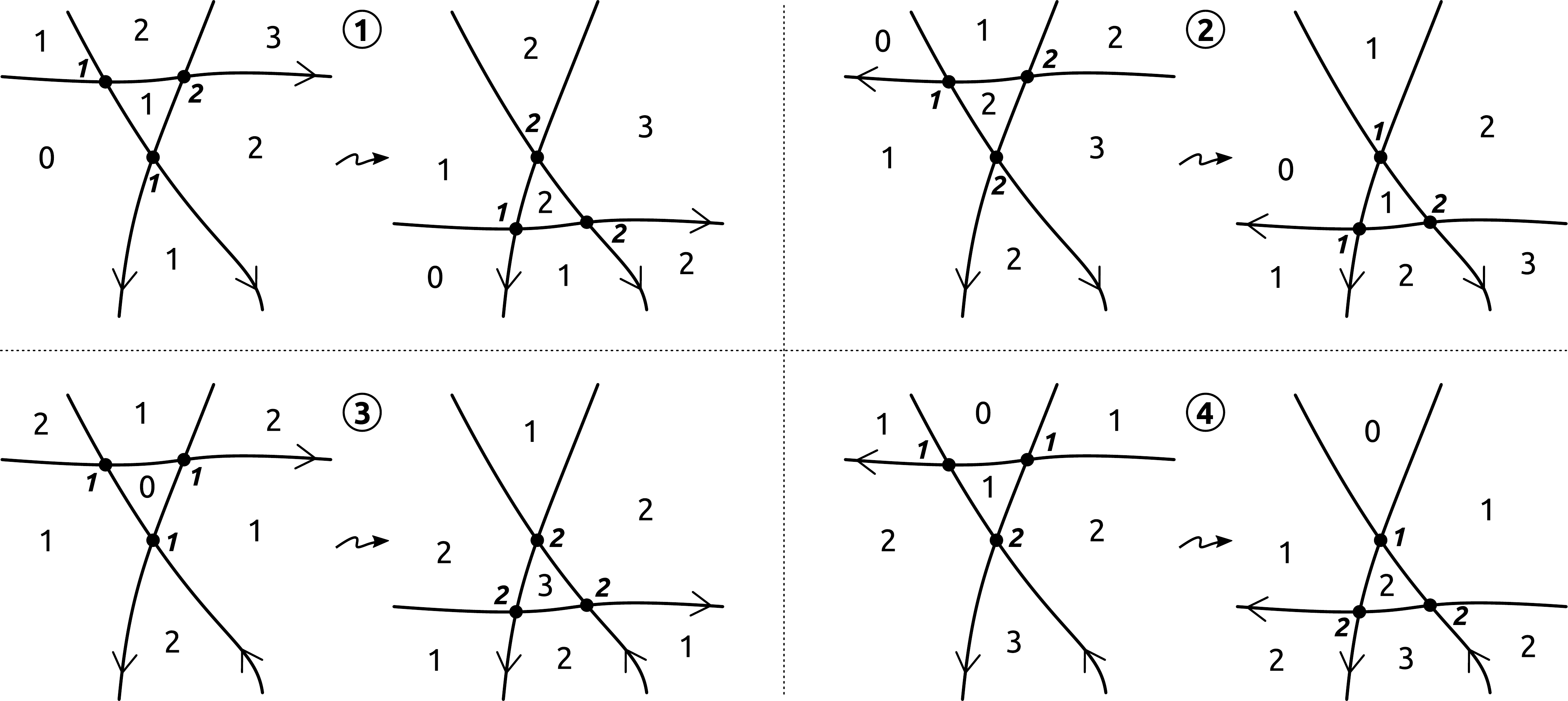}
		\caption{All four different cases of triple point crossings we need to consider. Labeled are the winding numbers for each connected component. The bold and cursive numbers are the indices of the double points.}
		\label{fig:roteqtrip}
	\end{figure}

	As before, we calculate the change when going from left to right in the pictures in Figure~\ref{fig:roteqtrip}:

	\begin{itemize}
		\item Case 1: $\underbrace{(2 - 1)}_{\text{change of} \sum \limits_{C \in \Gamma_K} \omega_C(K)} - \underbrace{((1 + 2 + 2) - (1 + 1 + 2))}_{\text{change of } \sum \limits_{p \in \mathcal{D}_K} \operatorname{ind}_p(K)} = 1 - (5 - 4) = 0$
		\item Case 2: $(1 - 2) - ((1 + 1 + 2) - (1 + 2 + 2)) = -1 - (4 - 5) = 0$
		\item Case 3: $(3 - 0) - ((2 + 2 + 2) - (1 + 1 + 1)) = 3 - (6 - 3) = 0$
		\item Case 4: $(2 - 1) - ((1 + 2 + 2) - (1 + 1 + 2)) = 1 - (5 - 4) = 0$
	\end{itemize}

	Now by the \emph{Whitney--Graustein Theorem} we know that $K$ is regularly homotopic to the standard curve~$K_j$ (see Definition~\ref{def:standardcurves}) with the same rotation number as $K,$ i.e.\@ to $K_{\operatorname{rot}(K)}.$ We know the equation is true for all standard curves and that it stays invariant under regular homotopies, which shows the proposition.
\end{proof}

\begin{remark}
	Earlier in Chapter~\ref{subsubsec:viro} we will introduce Viro's formula, a tool to calculate $J^+$ (yet to be introduced). It can be proven in \textnormal{exactly} the same way as the proof we just did. The full proof in this fashion can also be found on page 7 of \textnormal{The $J^{2+}$-Invariant for Pairs of Generic Immersions} by \textnormal{Hanna Häußler} \cite{hanna:paper}. Of course it is also proven in the cited literature, but differently. If you liked this proof, try proving Viro's formula later.
\end{remark}

Let us formulate a corollary from Lemma~\ref{lem:roteqsplit} that we will later need for the proof of Theorem~\ref{th:interiorsum} in Chapter~\ref{sec:intsum}.

\begin{thesiscorollary}
	\label{cor:dpscompsnumber}
	Let $K$ be any generic immersed loop -- i.e.\@ all self-intersections are transverse double points -- $n_K$ the number of double points of $K$ and $\Gamma_K$ its connected components. Then:
	\begin{equation*}
		| \Gamma_K | = n_K + 2
	\end{equation*}

	One of the connected components is unbounded and all the other are bounded.
\end{thesiscorollary}

\begin{proof}
	By the \emph{Jordan Curve Theorem} we know that the standard curve $K_1$ divides the plane $\mathbb{C}$ into one bounded component and one unbounded component, so
	$$| \Gamma_{K_1} | = 2 = 0 + 2 = n_{K_1} + 2 .$$
	The same is true for $K_{-1} .$

	The standard curve $K_0$ is a circle with an added outer loop. Adding a single loop increases the number of double points and the number of (bounded) connected components by $1. $ So
	$$| \Gamma_{K_0} | = 3 = 1 + 2 = n_{K_0} + 2 .$$

	For $|j| > 1$ we have a circle with $|j| - 1$ single interior loops, as we explored in the proof of Proposition~\ref{prop:rotequation}. So $n_{K_j} = |j| - 1$ and
	$$| \Gamma_{K_j} | = (|j| - 1) + 2 .$$

	Now again by the \emph{Whitney--Graustein Theorem} we know that the immersion~$K$ is regularly homotopic to the standard curve~$K_j$ (see Definition~\ref{def:standardcurves}) with the same rotation number as $K,$ i.e.\@ to $K_{\operatorname{rot}(K)}.$ We know the equation is true for all standard curves and that it stays invariant under regular homotopies (use Lemma~\ref{lem:roteqsplit}), which shows the corollary.
\end{proof}

\subsection{Proof of the theorem}

\begin{proof}[Proof of Theorem~\ref{th:interiorsum}: $J^+$ of interior sums]
	\label{pr:interiorsum}

	We use Viro's formula to prove the theorem. Let us just put the immersion $K'$ into $C$ and cross-connect the arcs $A$ and $A'$ and observe which new double points and winding numbers the resulting immersion $K^{\times}$ has compared to $K .$

	For the double points we get the ones that were already present in $K,$ all double points from $K'$ as well as a new one because of the cross-connection, so
	$$ n_{K^{\times}} = n_K + n_{K'} + 1 .$$

	For the connected components and their winding numbers we observe that the connected components we add to $K$ are all connected components from $K' ,$ but without the unbounded one (with winding number $0$). Let us denote $\Gamma_{K'}$ without the unbounded connected component as $\Gamma_{K'}^* .$ The winding number of these connected components in $K^{\times}$ is increased by $\omega_C(K) .$ Let us denote $b \vcentcolon= \omega_C(K)$ for readability, so:
	\begin{align*}
		\sum \limits_{\zeta \in \Gamma_{K^{\times}}} (\omega_\zeta(K^{\times}))^2
		&= \sum \limits_{\zeta \in \Gamma_{K}} (\omega_\zeta(K))^2 + \! \sum \limits_{\zeta \in \Gamma_{K'}^*} (\omega_\zeta(K') + b)^2 \\
		&= \sum \limits_{\zeta \in \Gamma_{K}} (\omega_\zeta(K))^2 + \! \sum \limits_{\zeta \in \Gamma_{K'}^*} (\omega_\zeta(K')^2 + 2 b \cdot \omega_\zeta(K') + b^2) \\
		&= \sum \limits_{\zeta \in \Gamma_{K}} (\omega_\zeta(K))^2 + \! \sum \limits_{\zeta \in \Gamma_{K'}} (\omega_\zeta(K'))^2 + \! \sum \limits_{\zeta \in \Gamma_{K'}} (2 b \cdot \omega_\zeta(K')) + \underbrace{|\Gamma_{K'}^*|}_{= n_{K'} + 1} \! \cdot \, b^2 \\
		&= \sum \limits_{\zeta \in \Gamma_{K}} (\omega_\zeta(K))^2 + \! \sum \limits_{\zeta \in \Gamma_{K'}} (\omega_\zeta(K'))^2 + 2 b \! \sum \limits_{\zeta \in \Gamma_{K'}} (\omega_\zeta(K')) + (n_{K'} + 1) \cdot b^2
	\end{align*}
	with Corollary~\ref{cor:dpscompsnumber} used in the third line.

	For the double points and their indices we see that the single double point that was created by the cross-connection has index $b \vcentcolon= \omega_C(K) .$ For the double points from $K'$ we see that their index is increased by $b$ as the winding numbers of all connected components were increased by $b ,$ so:
	\begin{align*}
		\sum \limits_{p \in \mathcal{D}_{K^{\times}}} (\operatorname{ind}_p(K^{\times}))^2
		&= \sum \limits_{p \in \mathcal{D}_{K}} (\operatorname{ind}_p(K))^2 + \! \sum \limits_{p \in \mathcal{D}_{K'}} (\operatorname{ind}_p(K') + b)^2 + b^2 \\
		&= \sum \limits_{p \in \mathcal{D}_{K}} (\operatorname{ind}_p(K))^2 + \! \sum \limits_{p \in \mathcal{D}_{K'}} (\operatorname{ind}_p(K')^2 + 2 b \cdot \operatorname{ind}_p(K') + b^2) + b^2 \\
		&= \sum \limits_{p \in \mathcal{D}_{K}} (\operatorname{ind}_p(K))^2 + \! \sum \limits_{p \in \mathcal{D}_{K'}} (\operatorname{ind}_p(K'))^2 + \!\! \sum \limits_{p \in \mathcal{D}_{K'}} (2 b \cdot \operatorname{ind}_p(K')) + \underbrace{|\mathcal{D}_{K'}|}_{= n_{K'}} \cdot \, b^2 + b^2 \\
		&= \sum \limits_{p \in \mathcal{D}_{K}} (\operatorname{ind}_p(K))^2 + \! \sum \limits_{p \in \mathcal{D}_{K'}} (\operatorname{ind}_p(K'))^2 + 2 b \! \sum \limits_{p \in \mathcal{D}_{K'}} (\operatorname{ind}_p(K')) + (n_{K'} + 1) \cdot b^2
	\end{align*}

	Now with the preparations done, let us use these observations and calculate $J^+(K^{\times}) .$ We use Viro's formula (see Lemma~\ref{lem:viro}) for the first line and to get $J^+(K)$ and $J^+(K')$ and directly after we use Equation~\ref{eq:rotequation} to get $\operatorname{rot}(K') :$
	\begin{align*}
		J^+(K^{\times})
		&= 1 + n_{K^{\times}} - \sum \limits_{\zeta \in \Gamma_{K^{\times}}} (\omega_\zeta(K^{\times}))^2 + \sum \limits_{p \in \mathcal{D}_{K^{\times}}} (\operatorname{ind}_p(K^{\times}))^2 \\
		&= 1 + (n_K + n_{K'} + 1) \\
			&\quad \, - ( \sum \limits_{\zeta \in \Gamma_{K}} (\omega_\zeta(K))^2 + \sum \limits_{\zeta \in \Gamma_{K'}} (\omega_\zeta(K'))^2 + 2 b \sum \limits_{\zeta \in \Gamma_{K'}} (\omega_\zeta(K')) + (n_{K'} + 1) \cdot b^2 ) \\
			&\quad \, + \sum \limits_{p \in \mathcal{D}_{K}} (\operatorname{ind}_p(K))^2 + \sum \limits_{p \in \mathcal{D}_{K'}} (\operatorname{ind}_p(K'))^2 + 2 b \sum \limits_{p \in \mathcal{D}_{K'}} (\operatorname{ind}_p(K')) + (n_{K'} + 1) \cdot b^2 \\
		&= \underbrace{1 + n_K - \sum \limits_{\zeta \in \Gamma_{K}} (\omega_\zeta(K))^2 + \sum \limits_{p \in \mathcal{D}_{K}} (\operatorname{ind}_p(K))^2}_{= J^+(K)} \\
			&\quad \, + \underbrace{1 + n_{K'} - \sum \limits_{\zeta \in \Gamma_{K'}} (\omega_\zeta(K'))^2 + \sum \limits_{p \in \mathcal{D}_{K'}} (\operatorname{ind}_p(K'))^2}_{= J^+(K')} \\
			&\quad \, - \left( 2 b \sum \limits_{\zeta \in \Gamma_{K'}} (\omega_\zeta(K')) + (n_{K'} + 1) \cdot b^2 \right) + \left( 2 b \sum \limits_{p \in \mathcal{D}_{K'}} (\operatorname{ind}_p(K')) + (n_{K'} + 1) \cdot b^2 \right) \\
		&= J^+(K) + J^+(K') - 2 b \underbrace{\left( \sum \limits_{\zeta \in \Gamma_{K'}} (\omega_\zeta(K')) - \sum \limits_{p \in \mathcal{D}_{K'}} (\operatorname{ind}_p(K')) \right)}_{= \operatorname{rot}(K')} \\
		&= J^+(K) + J^+(K') - 2 b \cdot \operatorname{rot}(K')
	\end{align*}

	Now we replace $b$ back with $\omega_C(K)$ and see that the theorem is proven.
\end{proof}

\begin{remark}
	Originally the theorem had a very different proof, the sketch for it follows:
	
	It is possible to pull a small part of the arc $A$ all the way outside of the immersion $K ,$ so that this small part of $A$ is adjacent to the unbounded component -- we can skip this first step if $C$ was already adjacent to the unbounded component. Pulling out this small part of $A$ changed $J^+$ of the immersion by some value $d \in 2 \mathbb{N}$ for all the new double points that were created pulling the small part out, whenever the self-tangency was direct, that will later be cancelled out. Now we take the connected sum of the two immersions, connecting the arc $A'$ of~$K'$ with the pulled out part of $A .$ Then we do a little trick which results in $K'$ being inside the small part of~$A'$ and changes $J^+$ by $\pm 2 \cdot \operatorname{rot}(K') ,$ see Figure~\ref{fig:intorprooftrick}.

	\begin{figure}[h!]
		\centering
		\includegraphics[scale=0.6]{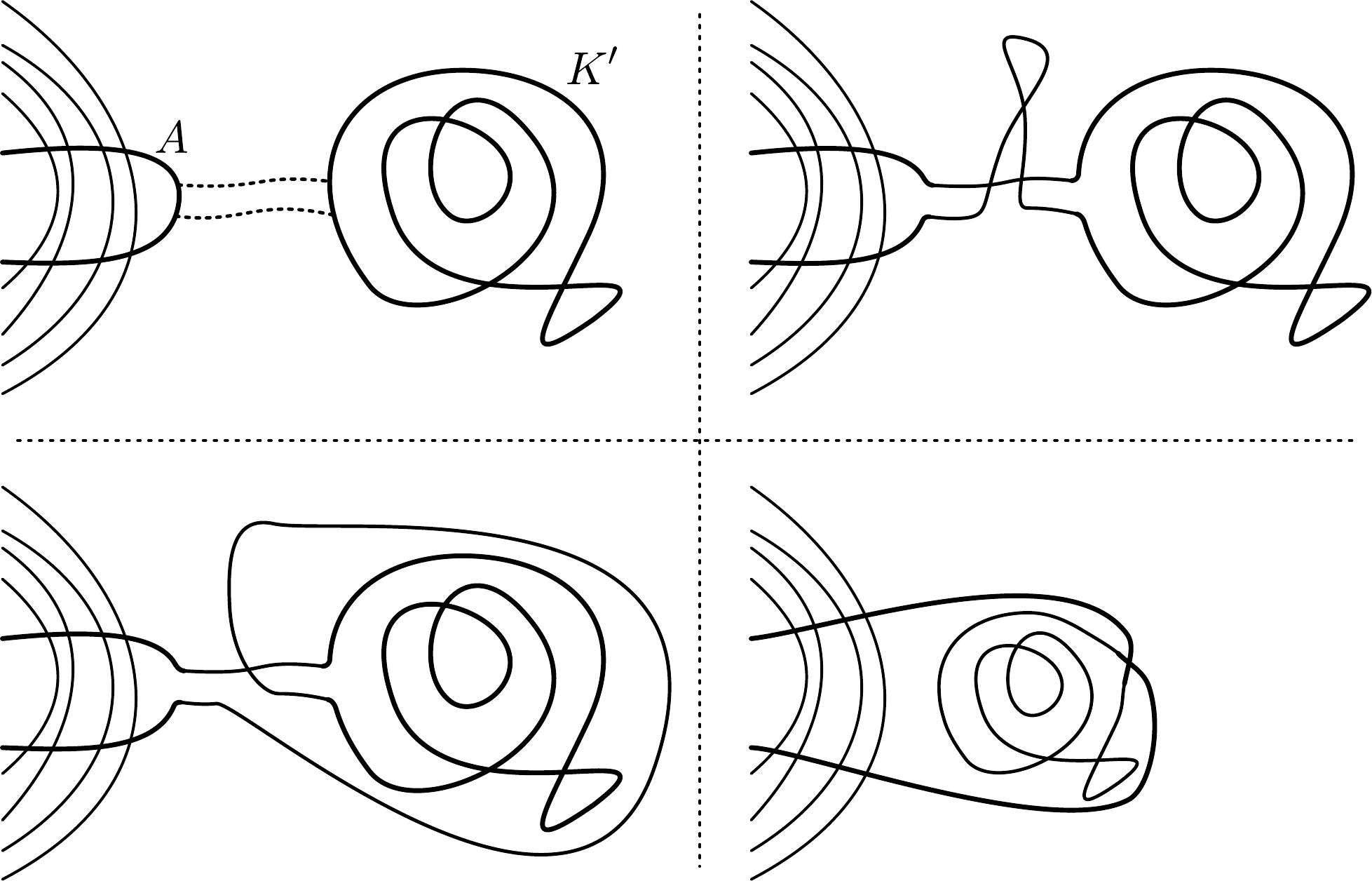}
		\caption{Illustration of a trick used in an earlier proof of Theorem~\ref{th:interiorsum}. The arcs at the left of the pictures hint at the arcs of $K$ that enclose the component $C$ and that have to be crossed by the small part of $A$ to reach into the unbounded component.}
		\label{fig:intorprooftrick}
	\end{figure}

	Then we pull the small part of $A$ all the way back to its original position as an arc of the component $C ,$ which also pulls the part that was originally the immersion $K' ,$ all the way into the component $C .$ This cancels out the change of $J^+$ by $d$ from before and also changes $J^+$ by $\pm 2 \cdot \omega_C(K) \cdot (\operatorname{rot}(K') \pm 1) .$ All instances of $\pm$ in this argument can be specifically determined, but need a lot more words. As does the whole proof.

	Formulating this whole idea properly turned out to be very tedious and temporarily involved an unproven assumption. Luckily, this original proof gave rise to the idea that the rotation number of the immersion~$K' ,$ which also played an important part in the original proof, should be equal to $\sum \limits_{C \in \Gamma_K} \omega_C(K) - \sum \limits_{p \in \mathcal{D}_K} \operatorname{ind}_p(K) .$ This then resulted in trying to prove Proposition~\ref{prop:rotequation}, which worked, and led to the realization that Theorem~\ref{th:interiorsum} can be proven using Viro's formula and this newly found equation for the rotation number. Let us all be happy that this lucky occurence saved us from what would have possibly been the worst proof of this paper.
\end{remark}

\begin{remark}
	Note that it does not matter which arc of the connected component~$C$ of the immersion~$K$ we cross-connect with the immersion~$K'$ as long as the orientation of $K'$ matches or is adjusted. We need to be careful with the sign of the rotation number of $K' .$ In the definition and in the theorem we assumed the orientation of $K'$ to already match. Taking another arc can require a change of orientation, which changes the sign of~$\operatorname{rot}(K')$ in the equation.
\end{remark}

\subsection{Corollaries}

Two noteworthy special cases of Theorem~\ref{th:interiorsum} are stated in the following two corollaries. The first one is also stated in \cite[Lemma 4 (with a different proof, not referring to a version of Theorem~\ref{th:interiorsum}, but instead using an argument around neighboring connected components)]{kai:paper}, but slightly different.

\begin{thesiscorollary}
	\label{cor:intloopkaiurs}
	Let $K$ be an arbitrary immersion.

	Consider any connected component $C$ of $\mathbb{C} \setminus K$ and any boundary arc $A$ of $C.$ Denote by~$C_{\text{adj}}$ the unique connected component of $\mathbb{C} \setminus K$ that also has $A$ as a boundary arc and $C_{\text{adj}} \neq C .$

	Denote the difference of winding numbers $\omega_C(K) - \omega_{C_{\text{adj}}}(K)$ as $\omega_{\text{adj}} \in \{ -1, 1 \} .$

	Denote as $K^{\text{o}}$ the immersion that we get if we add an interior loop (see Definition~\ref{def:intloops}) in the connected component~$C$ to the arc~$A$ of immersion~$K .$ Then
	$$J^+(K^{\text{o}}) = J^+(K) -2 \cdot \omega_C(K) \cdot \omega_{\text{adj}} .$$
\end{thesiscorollary}

\begin{proof}
	Adding an interior loop in $C$ to $A$ of the immersion $K$ is the same as the interior sum of the immersion~$K$ and a circle~$K_{\text{circle}}$ at boundary arc $A$ into the component $C ,$ where the orientation of the circle is adjusted to fit the cross-connection with the arc $A ,$ so $\operatorname{rot}(K_{\text{circle}}) = \omega_{\text{adj}}. $

	We use Theorem~\ref{th:interiorsum} and obtain
	\begin{align*}
		J^+(K^{\text{o}})
		&= J^+(K) + \underbrace{J^+(K_{\text{circle}})}_{= 0} -2 \cdot \omega_C(K) \cdot \underbrace{\operatorname{rot}(K_{\text{circle}})}_{= \omega_{\text{adj}}} \\
		&= J^+(K) -2 \cdot \omega_C(K) \cdot \omega_{\text{adj}} .
		\qedhere
	\end{align*}
\end{proof}

\vspace*{-0.6em}

\begin{thesiscorollary}
	\label{cor:mintloop}
	Let $K$ be an arbitrary immersion.

	Consider any connected component $C$ of $\mathbb{C} \setminus K$ and any boundary arc $A$ of $C.$ Denote by~$C_{\text{adj}}$ the unique connected component of $\mathbb{C} \setminus K$ that also has $A$ as a boundary arc and $C_{\text{adj}} \neq C .$

	Denote the difference of winding numbers $\omega_C(K) - \omega_{C_{\text{adj}}}(K)$ as $\omega_{\text{adj}} \in \{ -1, 1 \} .$

	Denote as $K^{\text{(m+1)}}$ the immersion that we get if we add an $(m+1)$-interior loop (see Definition~\ref{def:intloops}) in the connected component~$C$ to the arc~$A$ of immersion~$K .$ Then
	$$ J^+(K^{\text{(m+1)}}) = J^+(K) - (m + 1) (m + 2 \cdot \omega_C(K) \cdot \omega_{\text{adj}}) .$$
\end{thesiscorollary}

Before we prove this corollary, let us introduce the \emph{inner loop curves $A_j$} and calculate their $J^+$-value.

\begin{figure}[h!]
	\centering
	\includegraphics[scale=0.42]{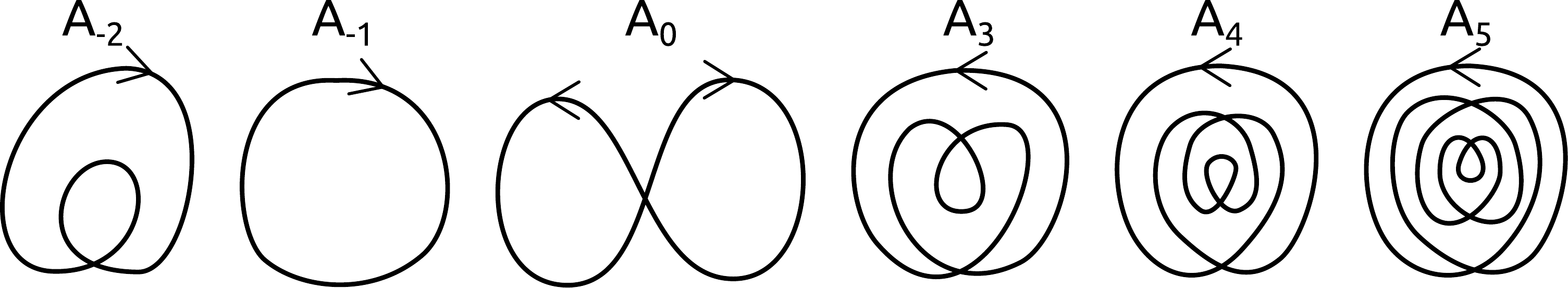}
	\caption{The inner loop curves. From left to right: $A_{-2}, A_{-1}, A_0, A_3, A_4, A_5 .$}
	\label{fig:innerloopcurves}
\end{figure}

\vspace*{-0.5em}

\begin{defi}[Inner loop curves $A_j$]
	\label{def:innerloopcurves}
	We call the immersions in Figure~\ref{fig:innerloopcurves} the \emph{inner loops curves} $A_j .$
	\begin{itemize}
		\item $A_j \vcentcolon= K_j$ for $j \in \{ -2, -1, 0, 1, 2 \}, $
		\item $\forall j, |j| > 1: A_j$ is a circle with an $(|j| - 1)$-interior loop and with $\operatorname{rot}(A_j) = j .$
	\end{itemize}
\end{defi}

In his paper \cite{viro:paper}, Viro proves Arnold's Conjecture. The conjecture states that for any arbitrary immersion~$K$ with $n$ many double points the lower bound $J^+(K) \ge -n^2 -n$ holds, and that for any number of double points $n$ the equation is attained by the inner loop curve~$A_{n+1} .$ We use our theorem for an alternate proof of the latter statement.

\begin{thesislemma}
	\label{lem:innerloopcurvesj}
	$$J^+(A_{n+1}) = -n^2 -n, \, \forall n \in \mathbb{N}$$
\end{thesislemma}

\begin{proof}
	We do a simple proof by induction using Corollary~\ref{cor:mintloop}.

	For the base case $n = 0$ we have:
	$$J^+(A_{n + 1}) = J^+(A_{0 + 1}) = J^+(A_1) = J^+(K_1) = 0 = (-1)^2 -1 \quad \checkmark$$

	For the induction step $n \to n + 1$ we get $A_{(n + 1) + 1} = A_{n + 2}$ from $A_{n + 1}$ if we add a single interior loop into the component $C$ with the highest winding number, which is $n + 1 .$ Using Corollary~\ref{cor:intloopkaiurs} we get:
	\begin{align*}
		J^+(A_{n + 2})
		&= \underbrace{J^+(A_{n + 1})}_{= -n^2 -n} -2 \cdot \underbrace{\omega_C(A_{n + 1})}_{= n + 1} \cdot \underbrace{\omega_{\text{adj}}}_{= 1} \\
		&= -n^2 -n - 2(n + 1) \\
		&= - ((n^2 + 2n + 1) + (n + 1)) \\
		&= - ((n+1)^2 + (n+1)) \\
		&= - (n+1)^2 - (n+1) \quad \checkmark
	\end{align*}

	Which proves the lemma.
\end{proof}

Now we can prove our corollary.

\begin{proof}[Proof of Corollary~\ref{cor:mintloop}]
	Adding an $(m+1)$-interior loop in $C$ to $A$ of the immersion $K$ is the same as the interior sum of the immersion~$K$ and the inner loop curve $A_{m+1}$ ($A_{-m-1}$ if $\omega_{\text{adj}} = -1$) at boundary arc $A$ into the component $C .$

	We use Theorem~\ref{th:interiorsum} and get the following. We write $A_{m+1}$ here in the first line, but the braces below the first line already account for the case if $\omega_{\text{adj}} = -1 ,$ where we would have to use $A_{-m-1}$:
	\begin{align*}
		J^+(K^{\text{(m+1)}})
		&= J^+(K) + \underbrace{J^+(A_{m+1})}_{= -m^2 - m} -2 \cdot \omega_C(K) \cdot \underbrace{\operatorname{rot}(A_{m+1})}_{= (m + 1) \cdot \omega_{\text{adj}}} \\
		&= J^+(K) -m^2 - m -2 \cdot \omega_C(K) \cdot (m + 1) \cdot \omega_{\text{adj}} \\
		&= J^+(K) - (m + 1) (m + 2 \cdot \omega_C(K) \cdot \omega_{\text{adj}}) .
		\qedhere
	\end{align*}
\end{proof}

The setup of the following corollary is a bit tedious, just like the definition of interior sums, but the statement itself is simple. It is rephrased informally after the corollary statement.

\begin{thesiscorollary}
	\label{cor:isolatedjchange}
	Let $K'$ and $K'_{\text{tr}}$ be two arbitrary immersions with the same rotation number, $\operatorname{rot}(K') = \operatorname{rot}(K'_{\text{tr}}) ,$ and each with at least one connected component $C'$ and $C'_{\text{tr}}$ adjacent to the unbounded connected component and with the same winding number~$\omega_{C'}(K') = \omega_{C'_{\text{tr}}}(K'_{\text{tr}}) \in \{ -1, 1 \} .$ Let $A'$ and $A'_{\text{tr}}$ be the boundary arcs of these connected components that are also boundary arcs of the unbounded component.

	Let $K$ be an arbitrary immersion.

	Consider any bounded connected component $C$ of $\mathbb{C} \setminus K$ and any boundary arc $A$ of $C.$

	Denote $K^{\times}$ as the interior sum of immersions $K, K'$ at boundary arcs $A, A'$ into the component $C$ of $\mathbb{C} \setminus K .$ And $K^{\times}_{\text{tr}}$ as the interior sum of immersions $K, K'_{\text{tr}}$ at boundary arcs $A, A'_{\text{tr}}$ into the component $C$ of $\mathbb{C} \setminus K .$

	Then:
	\begin{equation*}
		J^+(K^{\times}) = J^+(K^{\times}_{\text{tr}}) + (J^+(K') - J^+(K'_{\text{tr}})) .
	\end{equation*}
\end{thesiscorollary}
\vspace*{0.8em}

\begin{figure}[h!]
	\centering
	\includegraphics[scale=0.3]{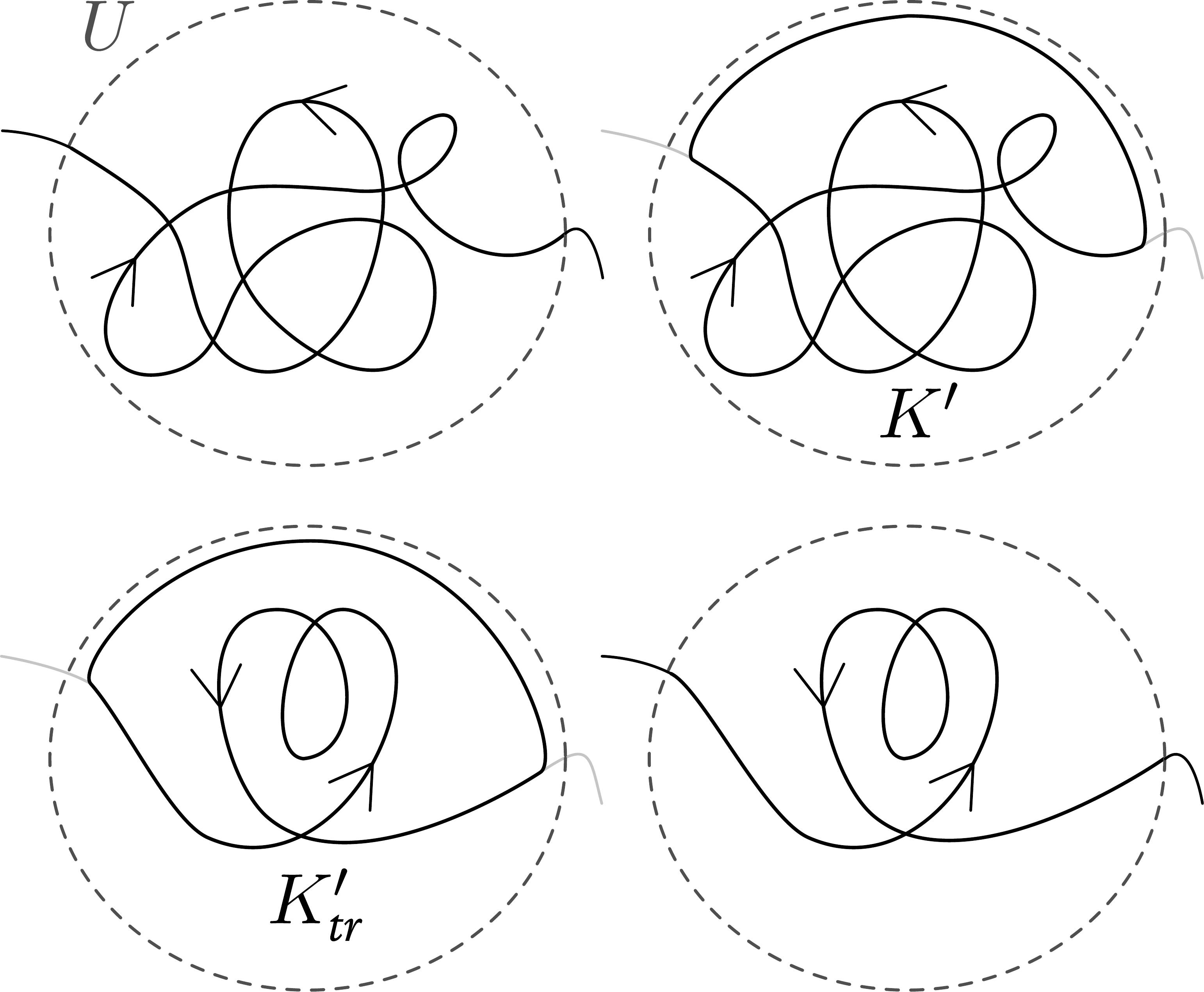}
	\caption{Illustration of Corollary~\ref{cor:isolatedjchange}. In the top left picture we select a simply connected area $U$ and construct $K'$ in the top right picture. We replace it by some other immersed loop with the same rotation number in the bottom left picture and reconnect it to the rest of the immersion in the bottom right picture. Here $\operatorname{rot}(K') = \operatorname{rot}(K'_{\text{tr}}) = 3, J^+(K') = -2, J^+(K'_{\text{tr}}) = -6.$}
	\label{fig:isolatedjchange}
\end{figure}

In other words this corollary states the following: Let $K^{\times} \subset \mathbb{C}$ be an arbitrary immersion and let $U \subset \mathbb{C}$ be a simply connected closed subset of $\mathbb{C}$ with $| \partial U \cap K^{\times} | = 2 ,$ so exactly two points of $K^{\times}$ are in the boundary of $U .$ Now ignore everything outside of $\overline{U} ,$ the closure of $U .$ Connect the two points of $K^{\times}$ in the boundary $\partial U$ of $U$ along $\partial U ,$ see Figure~\ref{fig:isolatedjchange}, and call this new immersion $K'$ that is the union of this connection (along $\partial U$) and the part of $K^{\times}$ in $U ,$ regularized (with smoothened edges).

Now apply any regular homotopies to $K' ,$ possibly with changes to its $J^+$-value, and call the new immersion $K'_{\text{tr}} .$ The only thing necessary is an outside arc of $K'_{\text{tr}}$ has to have the right orientation so that it can roughly be aligned with the part of $K'$ along $\partial U$ by applying a regular homotopy. Now remove that part along the border and rejoin the remaining part of $K'_{\text{tr}}$ with the rest of $K^{\times}$ and call this new immersion $K^{\times}_{\text{tr}} .$

Then we know that the difference of $J^+$ between~$K^{\times}_{\text{tr}}$ and $K^{\times}$ is equal to the difference of $J^+$ between $K'_{\text{tr}}$ and $K' .$

This guarantees that if two immersions are regularly homotopic -- i.e.\@ have the same rotation number -- and we can make out an isolated part of the immersion that needs to be changed to achieve the homotopy, then it is enough to analyze the change of $J^+$ in this isolated part of the immersion.

The way this is visualized and described in this remark might be easier to understand after the next corollary.

\begin{proof}
	With Theorem~\ref{th:interiorsum} and the following calculation the corollary follows:
	\begin{align*}
		J^+(K^{\times})
		&= J^+(K) + \underbrace{J^+(K')}_{= J^+(K') + J^+(K'_{\text{tr}}) - J^+(K'_{\text{tr}})} -2 \cdot \omega_C(K) \cdot \underbrace{\operatorname{rot}(K')}_{= \operatorname{rot}(K'_{\text{tr}})} \\
		&= J^+(K) + J^+(K'_{\text{tr}}) -2 \cdot \omega_C(K) \cdot \operatorname{rot}(K'_{\text{tr}}) + (J^+(K') - J^+(K'_{\text{tr}})) \\
		&= J^+(K^{\times}_{\text{tr}}) + (J^+(K') - J^+(K'_{\text{tr}}))
		\qedhere
	\end{align*}
\end{proof}

\begin{remark}
	The \emph{Whitney--Graustein Theorem} is also proven for isolated parts of an immersion -- i.e.\@ when a part of an immersion can be isolated with a simply connected area $U$ as in the corollary, then the isolated part is regularly homotopic to any immersion with the same rotation number via a regular homotopy that keeps the immersion outside of and at the border of $U$ unchanged during the homotopy -- which immediately proves Corollary~\ref{cor:isolatedjchange}.
\end{remark}

Now of course if we do not like to connect the immersions with a cross-connection, but instead with a tunnel, like at the connected sum, the following corollary gives us the formula for that.

\begin{figure}[h!]
	\centering
	\includegraphics[scale=0.35]{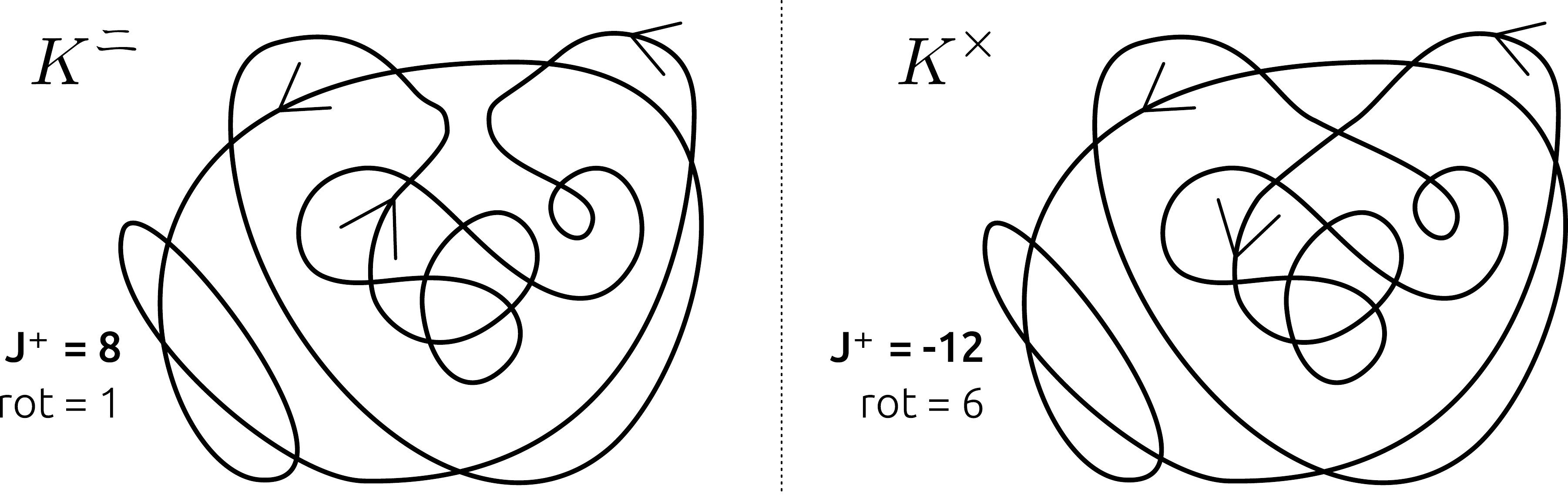}
	\caption{Tunnel-connected interior sum of the two immersions from Figure~\ref{fig:interiorsum}, pictured on the left. Note that this time the immersions are not connected with a cross-conection, but with two segments that do not intersect. Picture on the right shows the standard cross-connected interior sum for comparison.}
	\label{fig:interiortunnelsum}
\end{figure}

\begin{thesiscorollary}[$J^+$ of tunnel-connected interior sums]
	\label{cor:interiortunnelsum}
	Let $K, K', C, A, A', \omega_{\text{adj}}$ be the same as in Definition~\ref{def:interiorsum} and $K'$ already the correct orientation for the interior sum.

	Now denote by $K \jpniu$ the immersion that we construct like this: put $K'$ into $C$ without intersections between $K$ and $K'$ -- maybe it needs to be drawn smaller in visualizations -- and connect the arcs $A$ and $A'$ without intersections in the connection. To do this, the orientation of~$K'$ has to be changed. See Figure~\ref{fig:interiortunnelsum} for an illustration.

	Then
	\begin{equation*}
		J^+(K \jpniu) = J^+(K) + J^+(K') + 2 \cdot \omega_C(K) \cdot (\operatorname{rot}(K') - \omega_{\text{adj}}) .
	\end{equation*}
\end{thesiscorollary}

\begin{proof}
	Denote by $K^-$ the immersion $K'$ with opposite orientation.	First we take the connected sum of $K^-$ at the arc $A'$ and $K_0$ (figure eight) as illustrated in Figure~\ref{fig:interiortunnelsumpr}. We denote this new immersion as~$K^\infty ,$ its new boundary arc as $A^\infty$ and note that its rotation number is $\operatorname{rot}(K^\infty) = \operatorname{rot}(K^-) + \omega_{\text{adj}}$ and that $J^+(K^\infty) = J^+(K^-) = J^+(K') .$

	\begin{figure}[h!]
		\centering
		\includegraphics[scale=0.35]{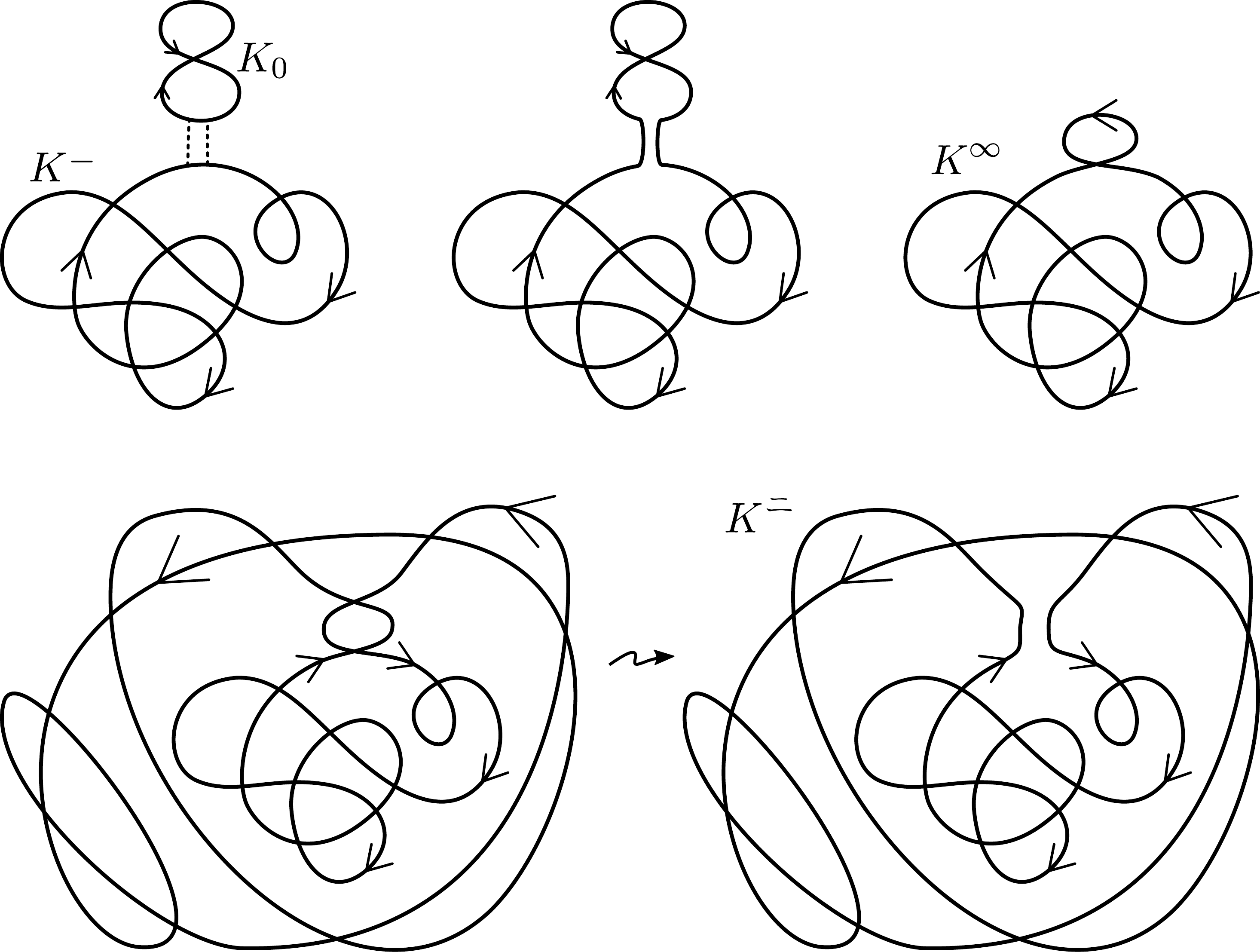}
		\caption{In the first row we take the connected sum of $K^-$ and $K_0$ at the arc $A'$ and apply a regular homotopy to get $K^\infty .$ The bottom left picture shows the interior sum of $K$ and $K^\infty$ and to get to the bottom right picture of $K \jpniu$ we apply a homotopy that goes through an inverse self-tangency.}
		\label{fig:interiortunnelsumpr}
	\end{figure}

	Now we just take the interior sum of immersions $K, K^\infty$ at boundary arcs $A, A^\infty$ into the component~$C$ of $\mathbb{C} \setminus K .$ In a moment we calculate what the $J^+$-value of the resulting immersion is, but first we do one last step to get $K \jpniu$ so that we can calculate its $J^+$-value.

	We apply a homotopy that -- through an inverse self-tangency (which does not change the value of $J^+$) -- removes the connected component that we added when we created $K^\infty ,$ see the last two pictures of Figure~\ref{fig:interiortunnelsumpr}, and end up with $K \jpniu .$

	With this and Theorem~\ref{th:interiorsum} we can calculate
	\begin{align*}
		J^+(K \jpniu)
		&= J^+(K) + \!\!\! \underbrace{J^+(K^\infty)}_{= J^+(K^-) \, = J^+(K')} \!\!\! -2 \cdot \omega_C(K) \,\, \cdot \!\!\!\!\! \underbrace{\operatorname{rot}(K^\infty)}_{= \underbrace{\operatorname{rot}(K^-)}_{= - \operatorname{rot}(K')} + \omega_{\text{adj}}} \\
		&= J^+(K) + J^+(K') - 2 \cdot \omega_C(K) \cdot ( - \operatorname{rot}(K') + \omega_{\text{adj}} ) \\
		&= J^+(K) + J^+(K') + 2 \cdot \omega_C(K) \cdot (\operatorname{rot}(K') - \omega_{\text{adj}})
		\qedhere
	\end{align*}
\end{proof}

\newpage

\section*{Appendix}
\addcontentsline{toc}{section}{\hspace{1.4em}Appendix}

\fancyhead[RO, LE]{Appendix}
\fancyhead[LO, RE]{}

Welcome to the appendix. The tex files and all image files of this paper can be found in the Git repository on \href{https://gitlab.com/CptMaister/paper-intro-to-j-plus}{\textbf{gitlab.com/CptMaister/paper-intro-to-j-plus}} and are available for anyone to use and modify.

\subsection*{Solutions}
\addcontentsline{toc}{subsection}{\hspace{2.3em}Solutions}

\subsubsection*{1. Ambiguous immersion (page \pageref{ex:01})}

\begin{figure}[h!]
	\centering
	\includegraphics[scale=0.42]{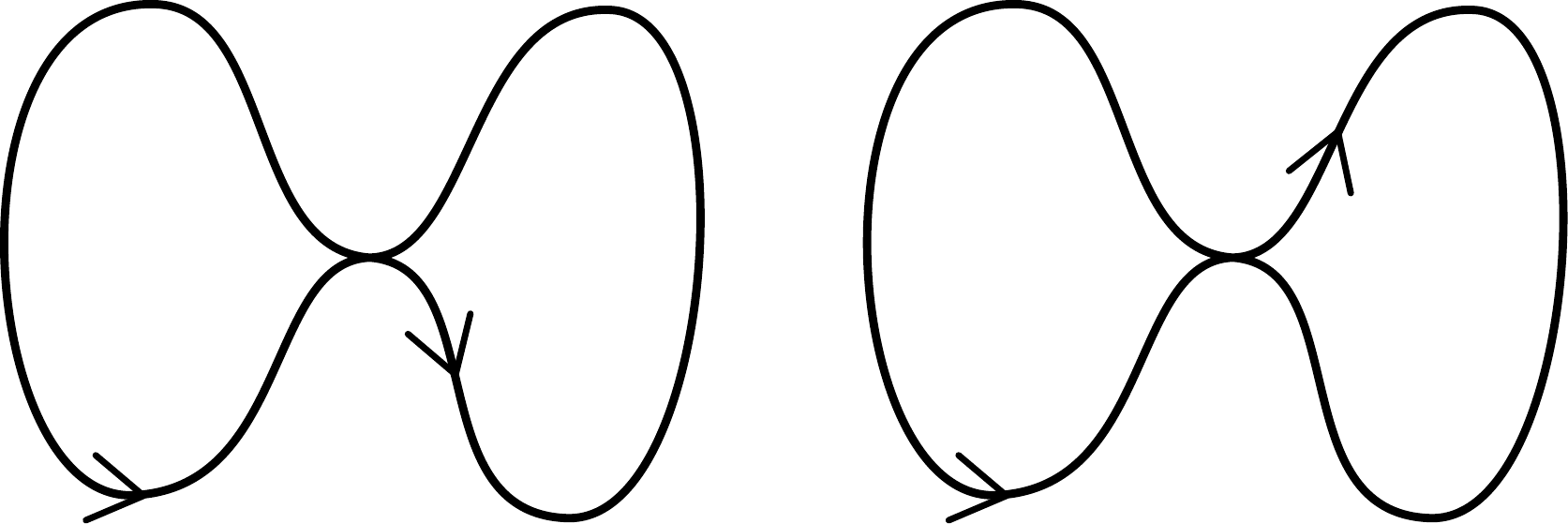}
	\caption{}
	\label{fig:sol01}
\end{figure}

Because of the tangential intersection of the immersion and the lack of an arrow indicating the orientation on the right half of the immersion, any of the two orientations of Figure~\ref{fig:sol01} are possible.

\subsubsection*{2. Not immersed loops (page \pageref{ex:02})}

The first, third, fourth and last curve is each not closed. The second one has edges, so it is not regular, i.e.\@ it is either not differentiable at the edges or the derivative vanishes (equals $0$) at the edges. The fifth picture is a mix of several not closed curves.

For any of these curves it can be argued that they are in fact closed and instead of being interrupted at some points, the curve just goes back the same way it came from in a sort of multiple cover. But then it is either not smooth or has vanishing derivative. The fourth picture (turtle) is then the only one that cannot be a single closed curve, as it is made up of at least three disjoint curves.

\subsubsection*{3. Fixing immersions (page \pageref{ex:03})}

\begin{figure}[h!]
	\centering
	\includegraphics[scale=0.42]{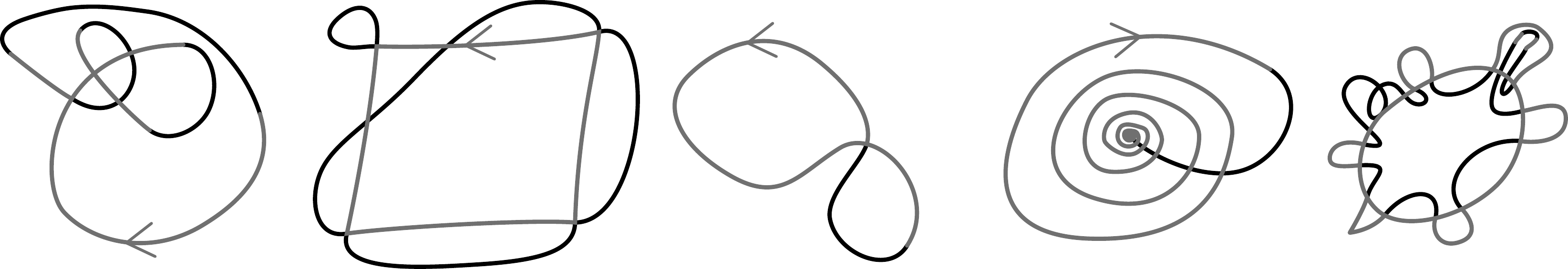}
	\caption{}
	\label{fig:sol03}
\end{figure}

Figure~\ref{fig:sol03} shows a suggestion for each of the curves, except for the last one. If this paper is printed, then the last curve can only be fixed by adding another paper next to the page and close the curve there.

\subsubsection*{4. Winding numbers and more (page \pageref{ex:04})}

To label the winding numbers and double point indices, we use the observations from Figure~\ref{fig:windingstep} and Remark~\ref{rem:indexdptwice}.

\begin{figure}[h!]
	\centering
	\includegraphics[scale=0.45]{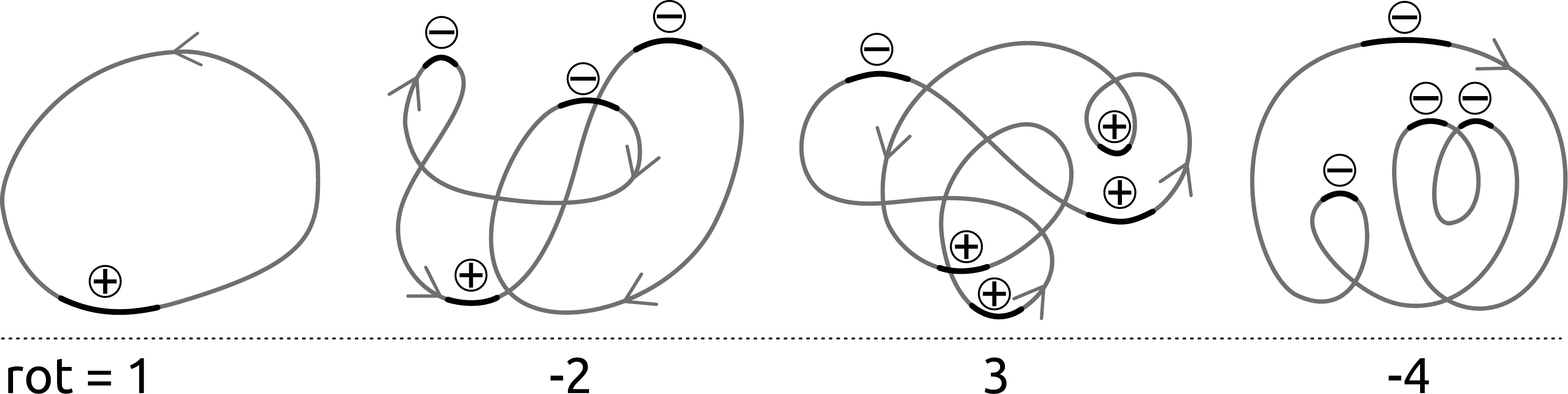}
	\caption{Happy and sad points for the immersions from Figure~\ref{fig:exrotwindind}. The difference of happy and sad points is equal to the rotation number.}
	\label{fig:sol04}
\end{figure}

For the rotation number there are several different ways to count it. One way is to mark all the points of the immersion where the tangent vector points to the right. Then let $a$ be the number of those points where the immersion curves up and $b$ the number of those points where the immersion curves down. Then the rotation number of the immersion is equal to $a - b .$ This method is sometimes interpreted as counting the smiles and frowns, or counting the happy and sad points of the immersion. In Figure~\ref{fig:sol04} all happy points are marked with a plus and all sad points with a minus.

\subsubsection*{5. Calculating $J^+$ basics 1 (page \pageref{ex:05})}

Let us denote the immersions of Figure~\ref{fig:jexamples} from left to right by $K^{5a}, K^{5b}, K^{5c}$ and $K^{5d} .$

We calculate $J^+$ of these immersions by applying regular homotopies to them until we reach a standard curve $K_j$ and keep track of the number of direct self-tangencies (abbreviate with \emph{dst} from here), as they change the value of $J^+ .$

\begin{figure}[h!]
	\centering
	\includegraphics[scale=0.4]{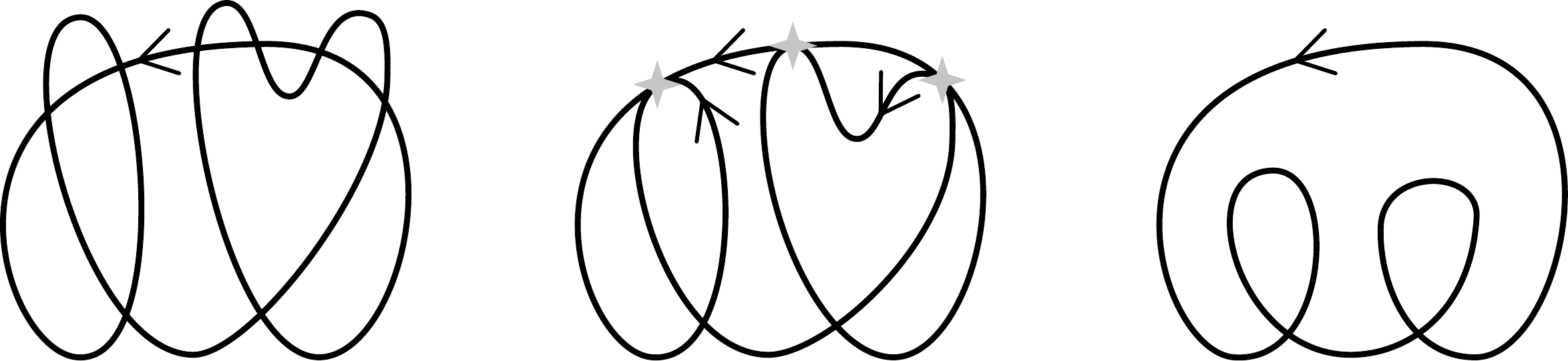}
	\caption{Homotopy from immersion~$K^{5a}$ to the standard curve $K_3 .$ The second picture shows the $3$ negative dst.}
	\label{fig:sol05a}
\end{figure}

For the first immersion $K^{5a}$, see Figure~\ref{fig:sol05a}, we arrive at $K_3$ with $0$ positive dst and $3$ negative dst. So the value of $J^+$ of the original immersion is
$$ J^+(K^{5a}) = J^+(K_3) - 2 (0 - 3) = -4 + 6 = 2 $$

\begin{figure}[h!]
	\centering
	\includegraphics[scale=0.4]{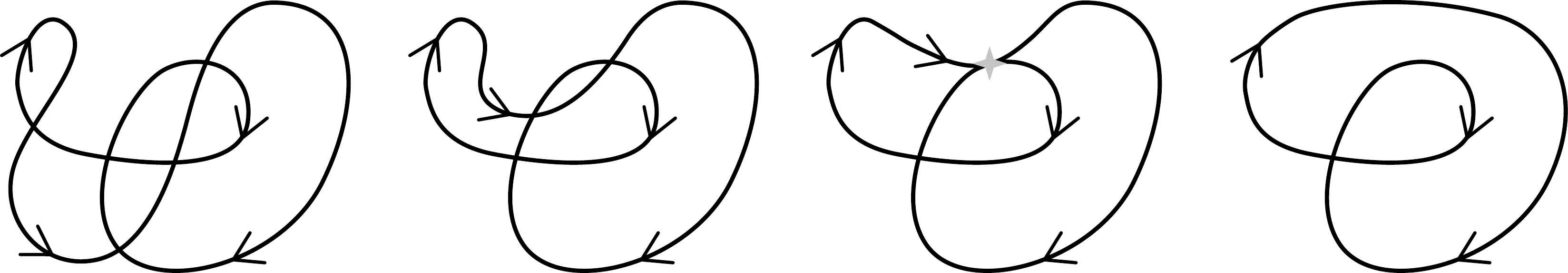}
	\caption{Homotopy from immersion~$K^{5b}$ to the standard curve $K_{-2} .$ There is one inverse self-tangency and one triple point crossing between the first and second picture. The third picture shows the negative dst.}
	\label{fig:sol05b}
\end{figure}

For the second immersion $K^{5b}$, see Figure~\ref{fig:sol05b}, we arrive at $K_{-2}$ with $0$ positive dst and $1$ negative dst. So the value of $J^+$ of the original immersion is
$$ J^+(K^{5b}) = J^+(K_{-2}) - 2 (0 - 1) = -2 + 2 = 0 $$

\begin{figure}[h!]
	\centering
	\includegraphics[scale=0.4]{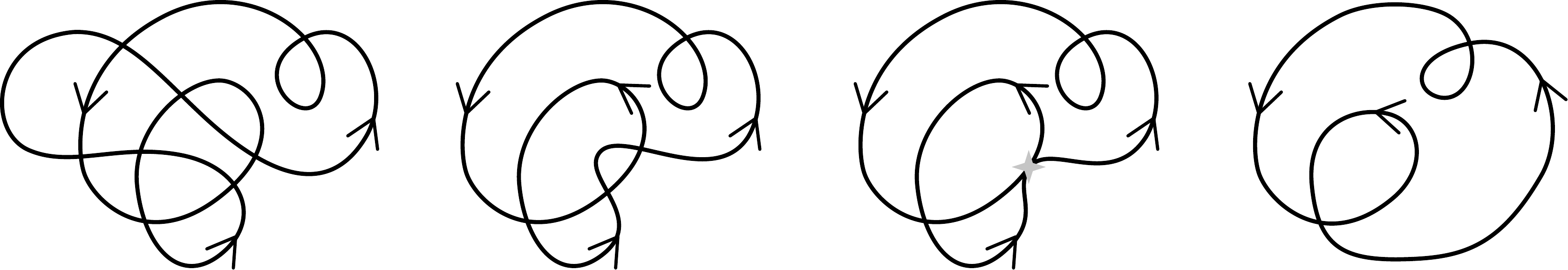}
	\caption{Homotopy from immersion~$K^{5c}$ to the standard curve $K_3 .$ There are two inverse self-tangencies between the first and second picture. The third picture shows the negative dst.}
	\label{fig:sol05c}
\end{figure}

For the third immersion $K^{5c}$, see Figure~\ref{fig:sol05c}, we arrive at $K_3$ with $0$ positive dst and $1$ negative dst. So the value of $J^+$ of the original immersion is
$$ J^+(K^{5c}) = J^+(K_3) - 2 (0 - 1) = -4 + 2 = -2 $$

\begin{figure}[h!]
	\centering
	\includegraphics[scale=0.4]{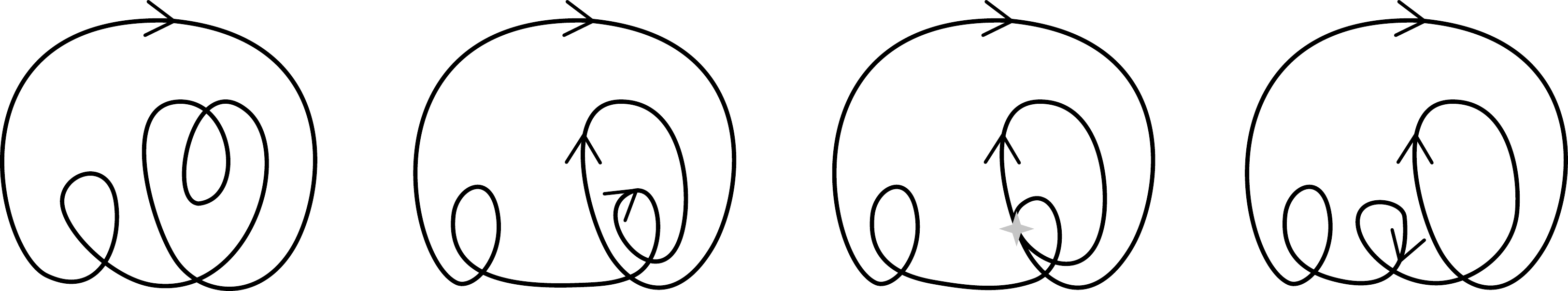}
	\caption{Homotopy from immersion~$K^{5d}$ to the standard curve $K_{-4} .$ The third picture shows the positive dst. There is one inverse self-tangency and one triple point crossing between the third and fourth picture.}
	\label{fig:sol05d}
\end{figure}

For the fourth immersion $K^{5d}$, see Figure~\ref{fig:sol05d}, we arrive at $K_{-4}$ with $1$ positive dst and $0$ negative dst. So the value of $J^+$ of the original immersion is
$$ J^+(K^{5d}) = J^+(K_{-4}) - 2 (1 - 0) = -6 - 2 = -8 $$

\subsubsection*{6. Calculating $J^+$ basics 2 (page \pageref{ex:06})}

Let us denote the immersions of Figure~\ref{fig:consum} by $K^{6a}$ (top left), $K^{6b}$ (top right) and $K^{6c}$ (bottom).

The first immersion $K^{6a}$ is the same as $K^{5d}$ from the previous exercise, see Figure~\ref{fig:sol05d}, with
$$J^+(K^{6a}) = J^+(K^{5d}) = -8 .$$

For the second immersion $K^{6b}$, we calculate $J^+$ by applying regular homotopies until we reach a standard curve $K_j$ and keep track of the number of direct self-tangencies (abbreviate with \emph{dst}), as they change the value of $J^+ .$

\begin{figure}[h!]
	\centering
	\includegraphics[scale=0.4]{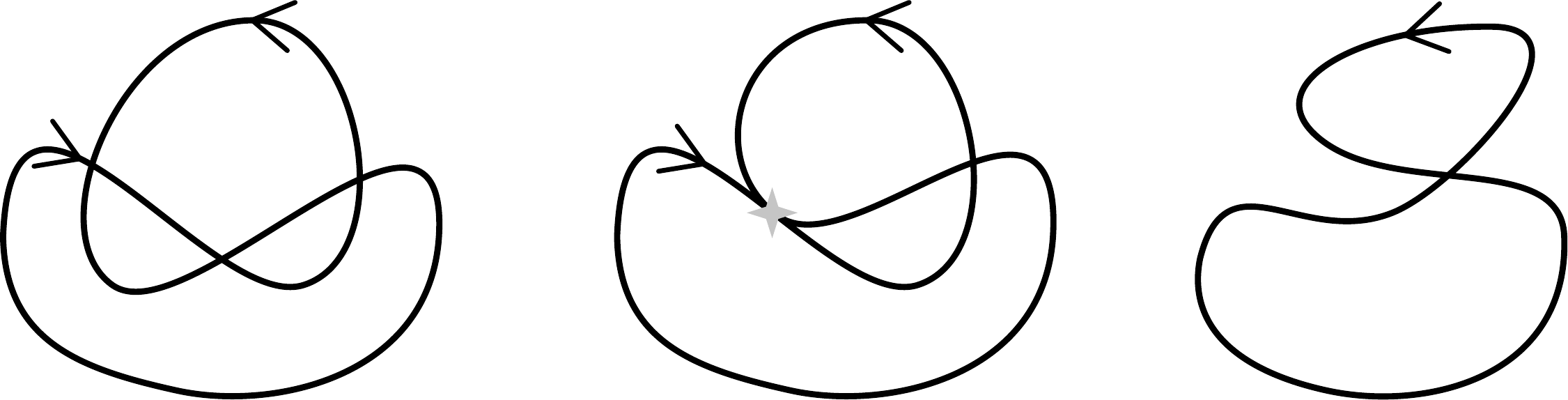}
	\caption{Homotopy from immersion~$K^{6b}$ to the standard curve $K_0 .$ The second picture shows the negative dst.}
	\label{fig:sol06b}
\end{figure}

We arrive at $K_0$ with $0$ positive dst and $1$ negative dst, see Figure~\ref{fig:sol06b}. So the value of $J^+$ of the original immersion is:
$$ J^+(K^{6b}) = J^+(K_0) - 2 (0 - 1) = 0 + 2 = 2 $$

For the bottom immersion $K^{6c}$ we can use the additivity of $J^+$ under connected sums:
$$ J^+(K^{6c}) = J^+(K^{6a}) + J^+(K^{6b}) = -8 + 2 = -6 $$

\subsubsection*{7. $K_j$ with intersecting interior loops (page \pageref{ex:07})}

\begin{figure}[h!]
	\centering
	\includegraphics[scale=0.4]{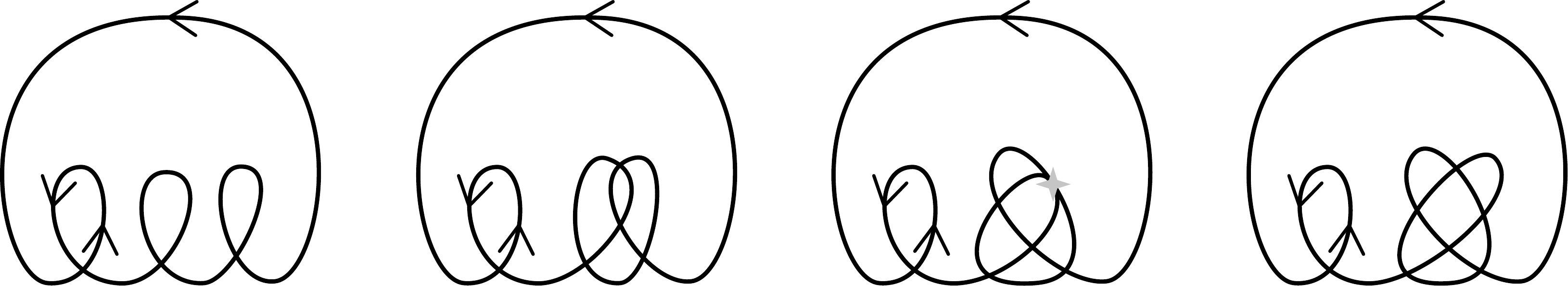}
	\caption{The standard curve $K_4$ in the first picture, with intersecting inner loops in the other pictures.}
	\label{fig:sol07}
\end{figure}

At first two interior loops intersect through an inverse self-tangency, see the second picture of Figure~\ref{fig:sol07}, so $J^+$ does not change.

Once they intersect, we can make them intersect through a positive direct self-tangency, see the third and fourth picture of Figure~\ref{fig:sol07}, which changes $J^+$ by $2 .$

\subsubsection*{8. $K_j$ with long interior loops (page \pageref{ex:08})}

\begin{figure}[h!]
	\centering
	\includegraphics[scale=0.4]{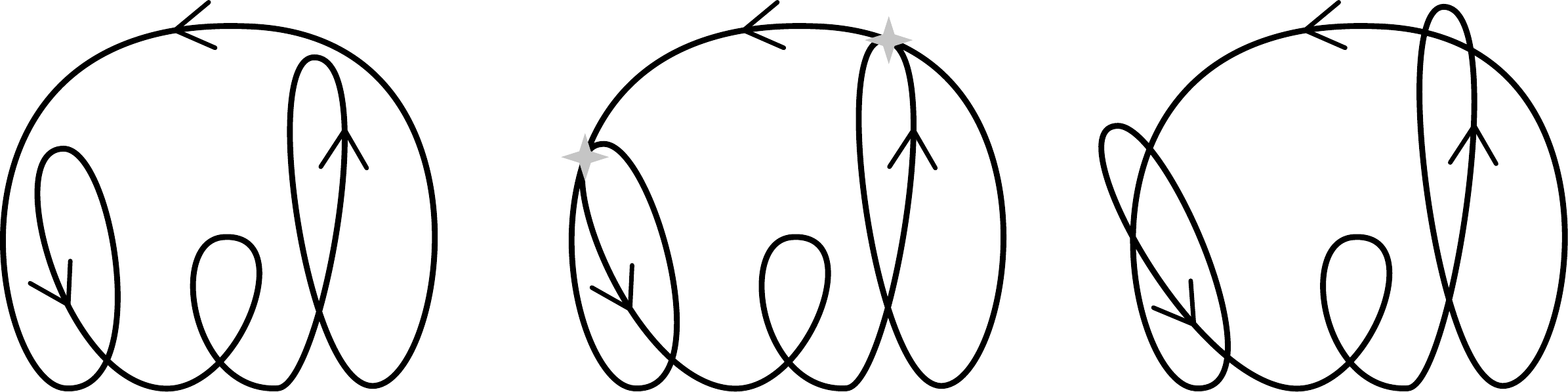}
	\caption{Standard curve $K_4$ with direct self-tangencies.}
	\label{fig:sol08}
\end{figure}

Any single intersection of a circle with its interior loops is through direct self-tangencies, which each changes $J^+$ by $2 .$ Figure~\ref{fig:sol08} shows the standard curve $K_4$ with two of its loops intersecting the upper part once each. Denote the immersion in the right picture as $K_4' ,$ then
$$J^+(K_4) = -6, \, J^+(K_4') = -2 .$$

\subsubsection*{9. Alternative calculation of $J^+(A_j)$ (page \pageref{ex:09})}

Let us denote a circle with a single interior loop by $A_2 ,$ a circle with a double interior loop by $A_3$ and a circle with a triple interior loop by $A_4 ,$ like in Definition~\ref{def:innerloopcurves}.

We know that $A_2 = K_2 ,$ so
$$J^+(A_2) = J^+(K_2) = -2 .$$

To calculate $J^+(A_3)$ we can start with $K_3$ and apply a regular homotopy that moves one of the two single interior loops into the other. During this homotopy there is $1$ negative direct self-tangency. So
$$J^+(A_3) = J^+(K_3) -2 (1) = -4 - 2 = -6 .$$

\begin{figure}[h!]
	\centering
	\includegraphics[scale=0.4]{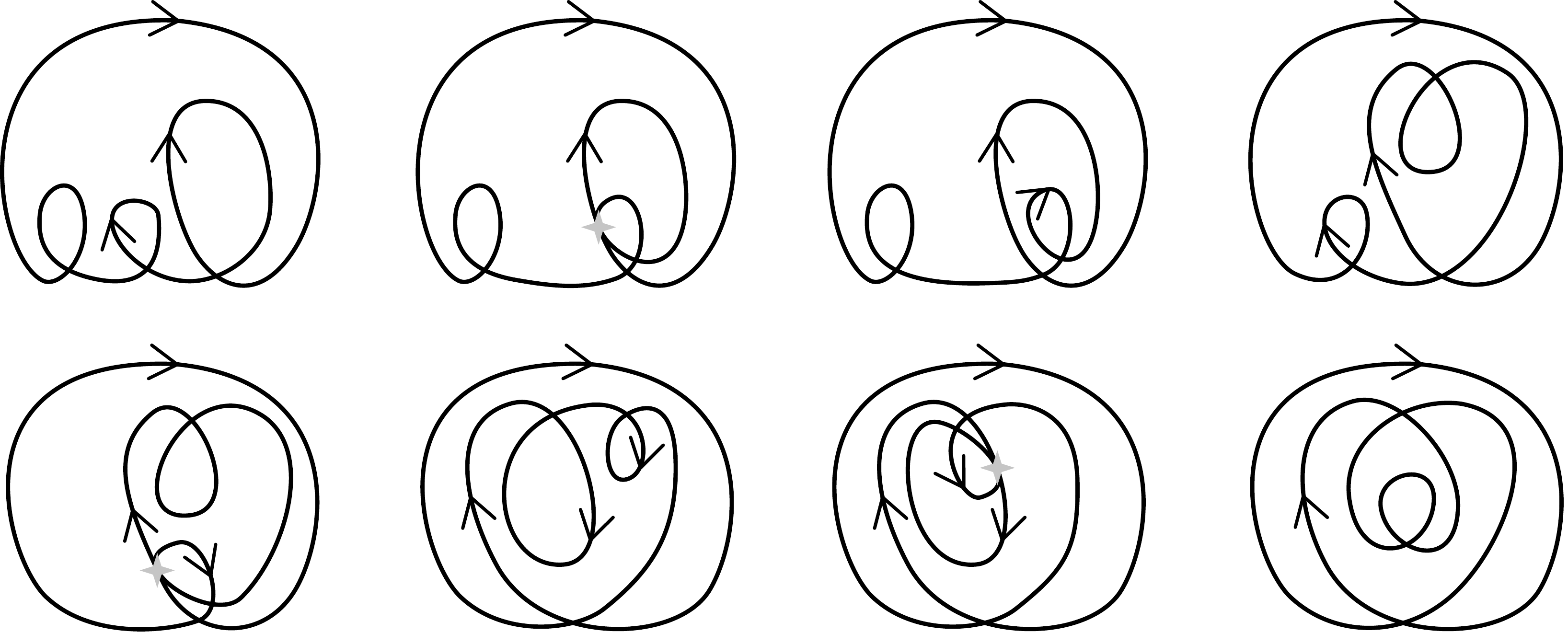}
	\caption{Homotopy from the standard curve~$K_4$ to the inner loop curve $A_4 .$}
	\label{fig:sol09b}
\end{figure}

And then we can calculate $J^+(A_4)$ similarly. We start with $K_4$ and apply a regular homotopy that moves one of the three single interior loops into the other, creating a double interior loop, with $1$ negative dst. Then we move the remaining single interior loop into the double interior loop, with $2$ negative dst, see Figure~\ref{fig:sol09b}. So
$$J^+(A_4) = J^+(K_4) -2 (1 + 2) = -6 - 6 = -12 .$$
	
This method can be used to calculate $J^+$ of any $A_j$ with $j > 1 ,$ which is a circle with a $(j-1)$-interior loop, and we can see that for any $j > 1 :$
\begin{align*}
	J^+(A_j)
	&= J^+(K_j) -2 \; \sum \limits_{k = 1}^{j-2} k \\
	&= -2 (j - 1) - 2 \; \frac{(j-2)(j-1)}{2} \\
	&= (j-1) (-2 - (j-2)) \\
	&= - j \ (j-1) \\
	&= - j^2 + j \\
	&= - (j-1)^2 - (j-1),
\end{align*}
which is identical to our result of Lemma~\ref{lem:innerloopcurvesj}.

\subsubsection*{10. $J^+$ and rotation number of interior sums (page \pageref{ex:10})}

First the rotation number of $K^{\times}$ and $K \jpniu .$ Let us use the denotions from Corollary~\ref{cor:interiortunnelsum} [$J^+$ of tunnel-connected interior sums] and its proof. We can calculate
$$\operatorname{rot}(K) = \operatorname{rot}(K') = 3, \quad \operatorname{rot}(K^{\infty}) = -2$$

With the happy-sad-points-approach from Figure~\ref{fig:sol04} it is easy to see that the rotation number of the cross-connected interior sum of two immersions is the sum of their rotation numbers. This is because the cross-connection does not remove or add any happy or sad points. It follows, that
$$\operatorname{rot}(K^{\times}) = \operatorname{rot}(K) + \operatorname{rot}(K') = 3 + 3 = 6, \quad
\operatorname{rot}(K \jpniu) = \operatorname{rot}(K) + \operatorname{rot}(K^{\infty}) = 3 - 2 = 1.$$

For the $J^+$-value of $K^{\times}$ and $K \jpniu$ we use the results of Theorem~\ref{th:interiorsum} and Corollary~\ref{cor:interiortunnelsum}. We know, that
$$J^+(K) = 2, \quad J^+(K') = -2 .$$

And with that, we can calculate
\begin{align*}
	J^+(K^{\times})
	&= J^+(K) + J^+(K') - 2 \cdot \omega_C(K) \cdot \operatorname{rot}(K') \\
	&= 2 - 2 - 2 \cdot 2 \cdot 3 \\
	&= -12 ,
\end{align*}
with $C$ the connected component of $K$ in Figure~\ref{fig:interiorsum}, and
\begin{align*}
	J^+(K \jpniu)
	&= J^+(K) + J^+(K') + 2 \cdot \omega_C(K) \cdot (\operatorname{rot}(K') - \omega_{\text{adj}}) \\
	&= 2 - 2 + 2 \cdot 2 \cdot (3 - 1) \\
	&= 8 .
\end{align*}

\subsubsection*{11. Calculating $J^+$ advanced (page \pageref{ex:11})}

Let us denote the immersions of Figure~\ref{fig:manyjexamples1} from left to right by~$K^a, K^b, K^c$ and $K^d ,$ the immersions of Figure~\ref{fig:manyjexamples2} from left to right by~$K^e, K^f, K^g$ and $K^h ,$ the immersions of Figure~\ref{fig:manyjexamples3} from left to right by~$K^p, K^q, K^r$ and $K^s ,$ and the immersions of Figure~\ref{fig:manyjexamples4} from left to right by~$K^t, K^u$ and $K^v .$

The immersion $K^a$ is regularly homotopic to $K_9$ without any self-tangencies at all, so
$$J^+(K^a) = J^+(K_9) = -16 .$$

\begin{figure}[h!]
	\centering
	\includegraphics[scale=0.4]{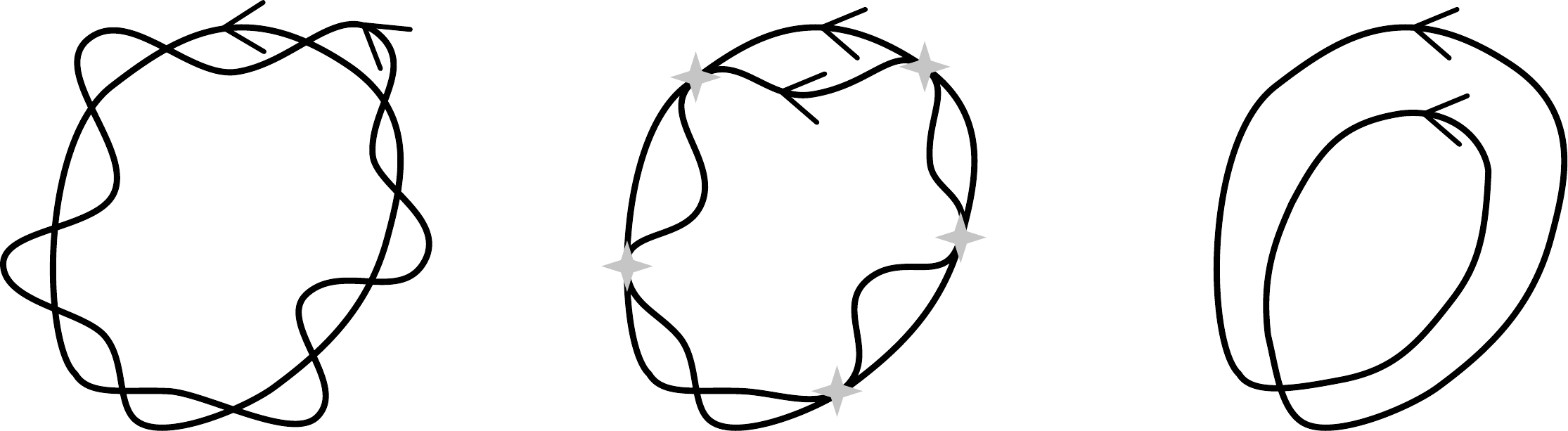}
	\caption{}
	\label{fig:sol11b}
\end{figure}

Immersion $K^b$ is a circle with a single interior loop that intersects the circle $5$ times through positive direct self-tangencies (abbreviate with \emph{dst}), see Figure~\ref{fig:sol11b}, so
$$J^+(K^b) = J^+(K_2) + 2 \cdot 5 = -2 + 10 = 8 .$$

\begin{figure}[h!]
	\centering
	\includegraphics[scale=0.4]{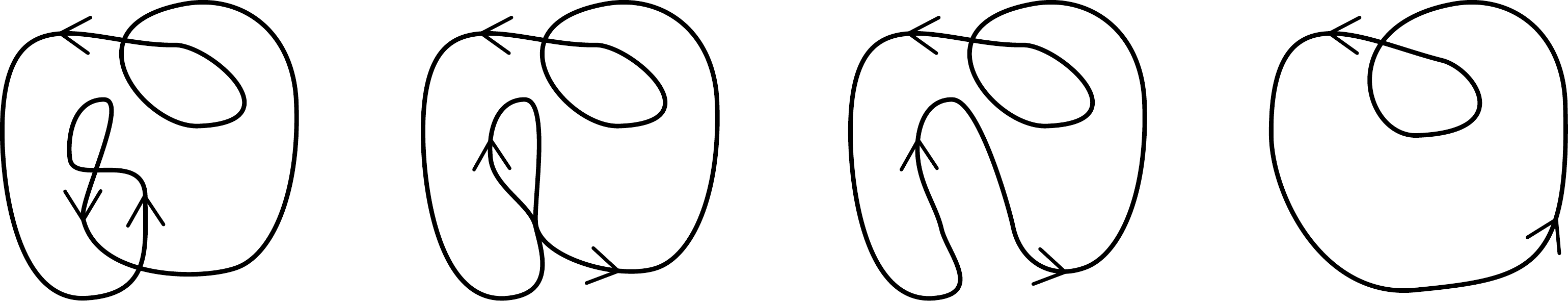}
	\caption{}
	\label{fig:sol11c}
\end{figure}

Immersion $K^c$ is regularly homotopic to $K_2$ through one inverse self-tangency, see Figure~\ref{fig:sol11c}, so
$$J^+(K^c) = J^+(K_2) = -2 .$$

\begin{figure}[h!]
	\centering
	\includegraphics[scale=0.4]{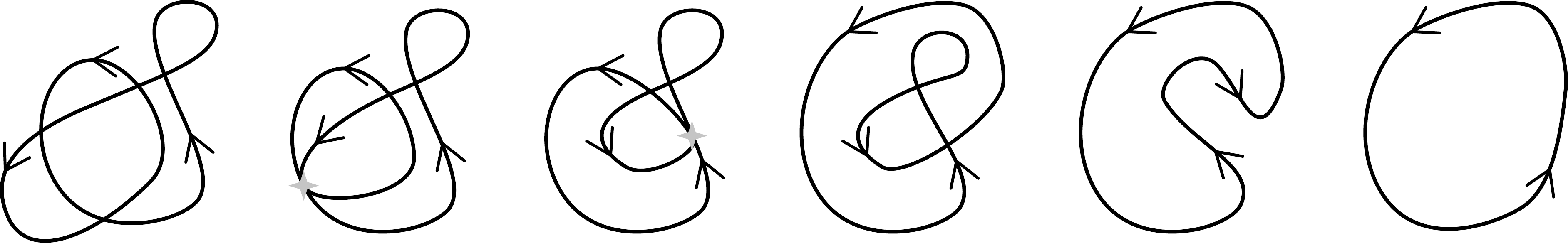}
	\caption{}
	\label{fig:sol11d}
\end{figure}

Immersion $K^d$ is regularly homotopic to $K_1 ,$ a circle, through $1$ negative dst and $1$ positive dst, see Figure~\ref{fig:sol11d}, so
$$J^+(K^d) = J^+(K_1) - 2(1 - 1) = 0 .$$

\begin{figure}[h!]
	\centering
	\includegraphics[scale=0.4]{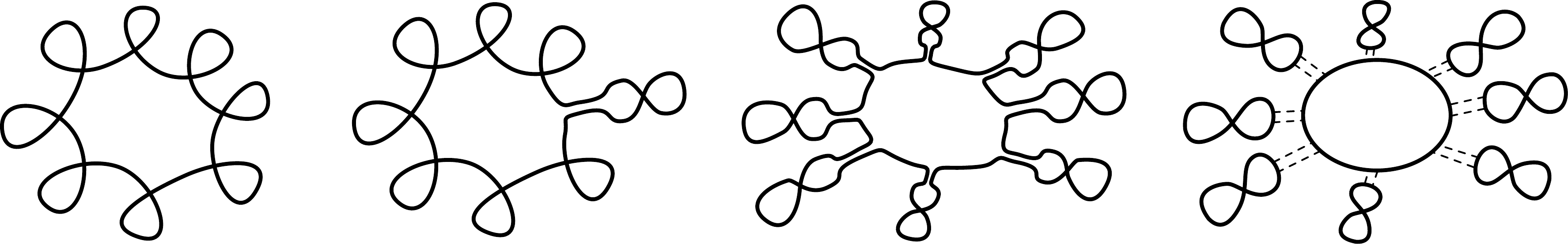}
	\caption{}
	\label{fig:sol11e}
\end{figure}

Immersion $K^e$ can be created by taking the connected sum of a circle, $K_1 ,$ and many figure eights,~$K_0 ,$ see Figure~\ref{fig:sol11e}, so by the connectivity of $J^+$ under connected sums we get
$$J^+(K^e) = J^+(K_1) + 8 \cdot J^+(K_0) = 0 + 8 \cdot 0 = 0 .$$

\begin{figure}[h!]
	\centering
	\includegraphics[scale=0.4]{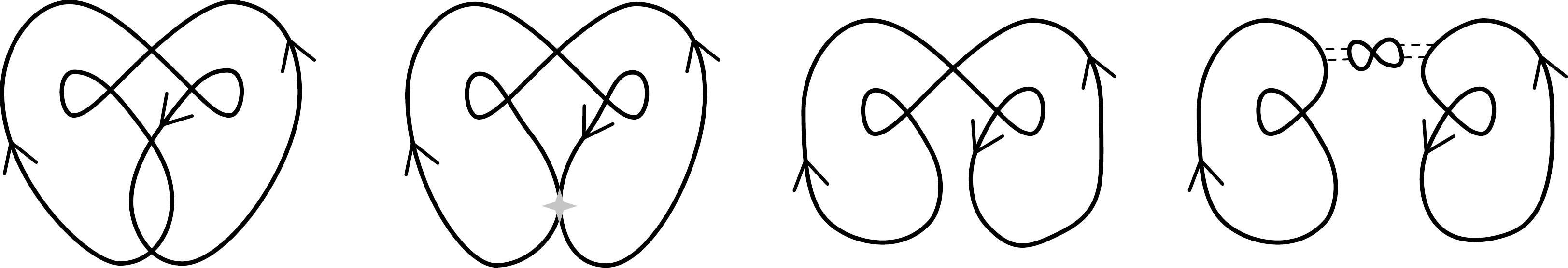}
	\caption{}
	\label{fig:sol11f}
\end{figure}

Immersion $K^f$ can be created by taking the connected sum of $K_{-2}, K_0$ and $K_2$ and then a regular homotopy with $1$ positive dst, see Figure~\ref{fig:sol11f} (from right to left), so by the connectivity of $J^+$ under connected sums we get
$$J^+(K^e) = J^+(K_{-2}) + J^+(K_0) + J^+(K_2) + 2 = -2 + 0 - 2 + 2 = -2 .$$

\begin{figure}[h!]
	\centering
	\includegraphics[scale=0.4]{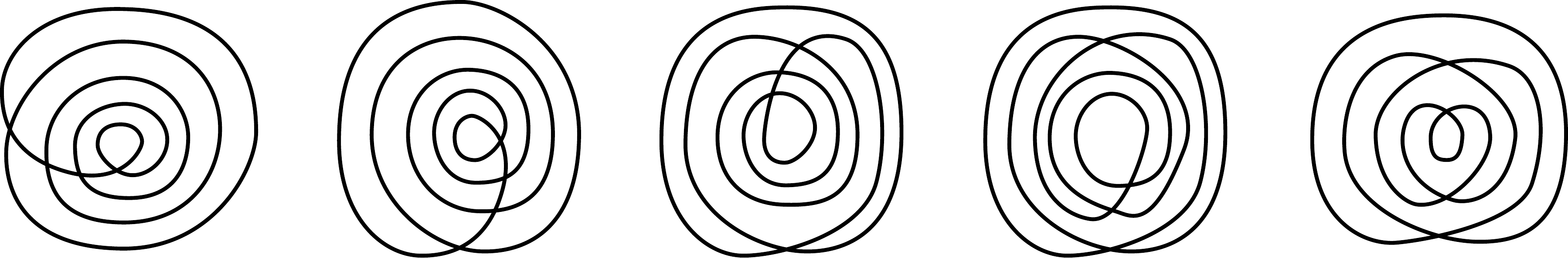}
	\caption{}
	\label{fig:sol11g}
\end{figure}

Immersion $K^g$ is regularly homotopic to $A_5 ,$ a circle with a $4$-interior loop, without any self-tangencies, see Figure~\ref{fig:sol11g}, so
$$J^+(K^g) = J^+(A_5) = - 5^2 + 5 = -20 .$$

\begin{figure}[h!]
	\centering
	\includegraphics[scale=0.35]{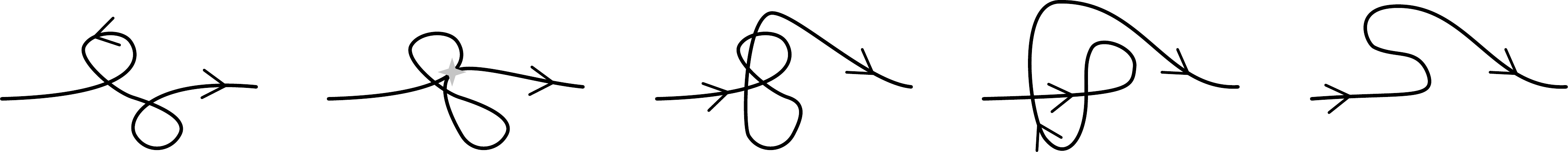}
	\caption{}
	\label{fig:sol11htrick}
\end{figure}

For immersion $K^h$ we need a little trick, which is illustrated in Figure~\ref{fig:sol11htrick}, that removes two counteroriented loops that are next to each other and increases $J^+$ by $2$ as there is one positive dst. 

\begin{figure}[h!]
	\centering
	\includegraphics[scale=0.4]{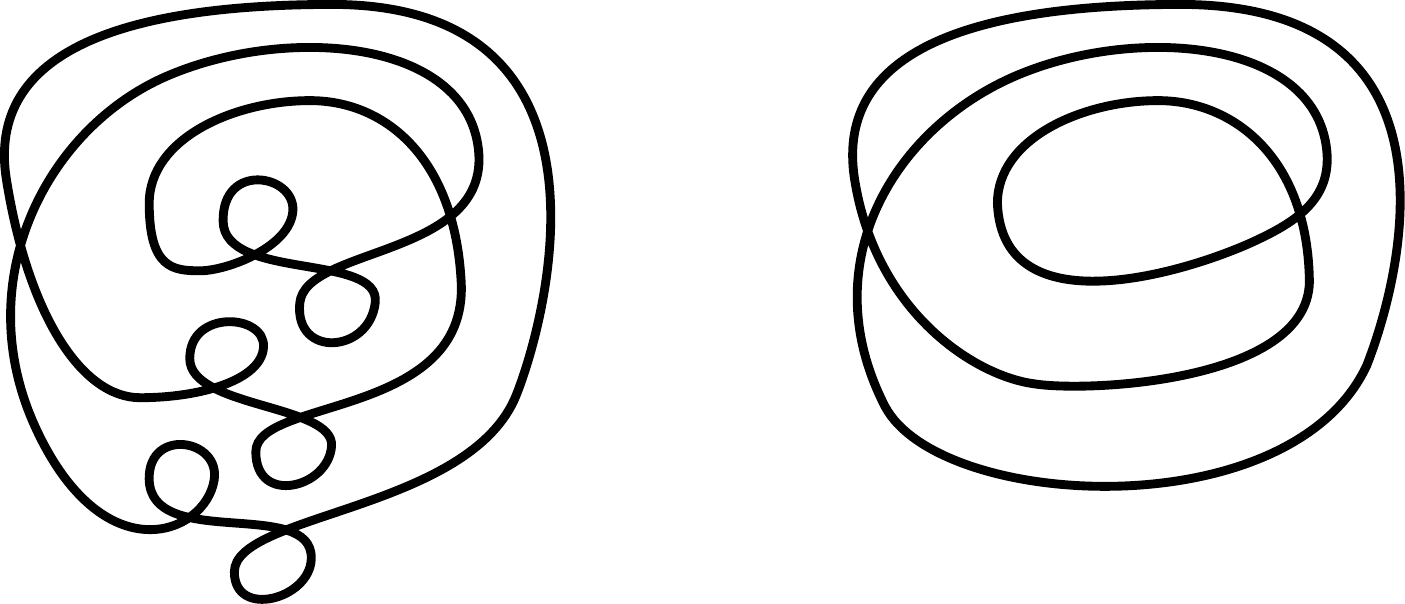}
	\caption{}
	\label{fig:sol11h}
\end{figure}

With this, we see that immersion $K^h$ is regularly homotopic to $A_3 ,$ a circle with a $2$-interior loop, after using the trick to remove counteroriented loops three times, see Figure~\ref{fig:sol11h}, so
$$J^+(K^h) = J^+(A_3) - 2 \cdot 3 = - 3^2 + 3 - 6  = -12.$$

\begin{figure}[h!]
	\centering
	\includegraphics[scale=0.4]{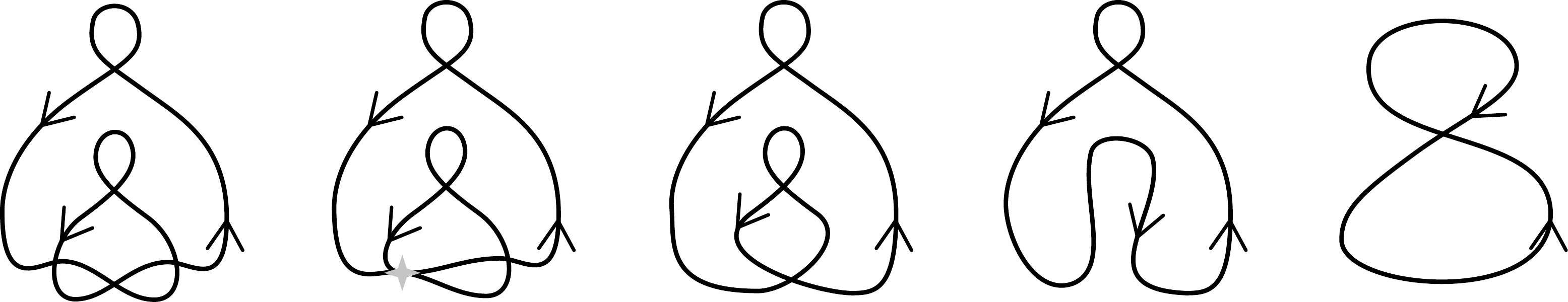}
	\caption{}
	\label{fig:sol11p}
\end{figure}

Immersion $K^p$ is regularly homotopic to the standard curve $K_0 ,$ the figure eight, through $1$ negative dst, see Figure~\ref{fig:sol11p}, so
$$J^+(K^p) = J^+(K_0) + 2 = 2 .$$

\begin{figure}[h!]
	\centering
	\includegraphics[scale=0.4]{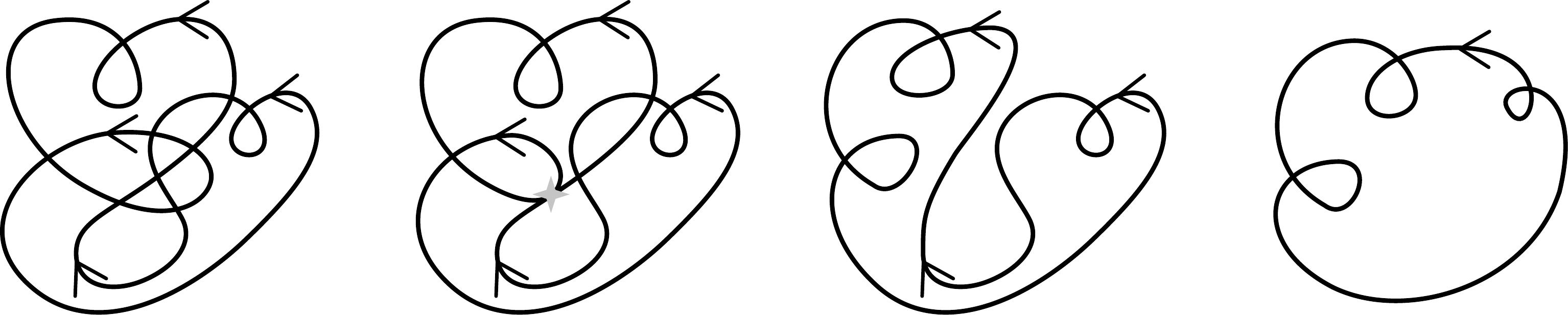}
	\caption{}
	\label{fig:sol11q}
\end{figure}

Immersion $K^q$ is regularly homotopic to the standard curve $K_4 ,$ a circle with three single interior loops, through $1$ negative dst, see Figure~\ref{fig:sol11q}, so
$$J^+(K^q) = J^+(K_4) + 2 = -6 + 2 = -4 .$$

\begin{figure}[h!]
	\centering
	\includegraphics[scale=0.36]{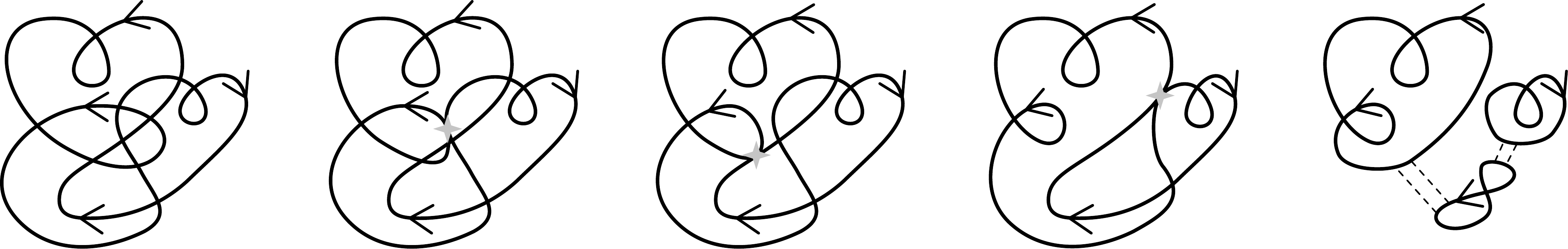}
	\caption{}
	\label{fig:sol11r}
\end{figure}

Immersion $K^r$ is regularly homotopic to the connected sum of the standard curves $K_3, K_0$ and $K_{-2} ,$ through $3$ negative dst, see Figure~\ref{fig:sol11r}, so
$$J^+(K^r) = J^+(K_3) + J^+(K_0) + J^+(K_{-2}) + 2 \cdot 3 = -4 + 0 - 2 + 6 = 0 .$$

Immersion $K^s$ is the connected sum of many figure eights, similar to the immersion $K^e, $ so
$$J^+(K^s) = 0 .$$

\begin{figure}[h!]
	\centering
	\includegraphics[scale=0.4]{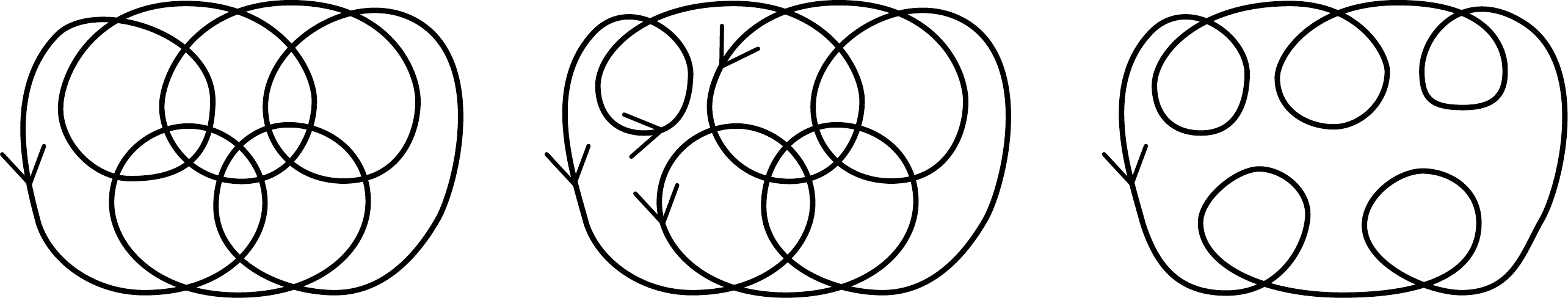}
	\caption{}
	\label{fig:sol11t}
\end{figure}

Immersion $K^t$ is regularly homotopic to the standard curve $K_6 ,$ a circle with five single interior loops, through only inverse self-tangencies, no direct self-tangencies, see Figure~\ref{fig:sol11t}, so
$$J^+(K^t) = J^+(K_6) = -10 .$$

\begin{figure}[h!]
	\centering
	\includegraphics[scale=0.4]{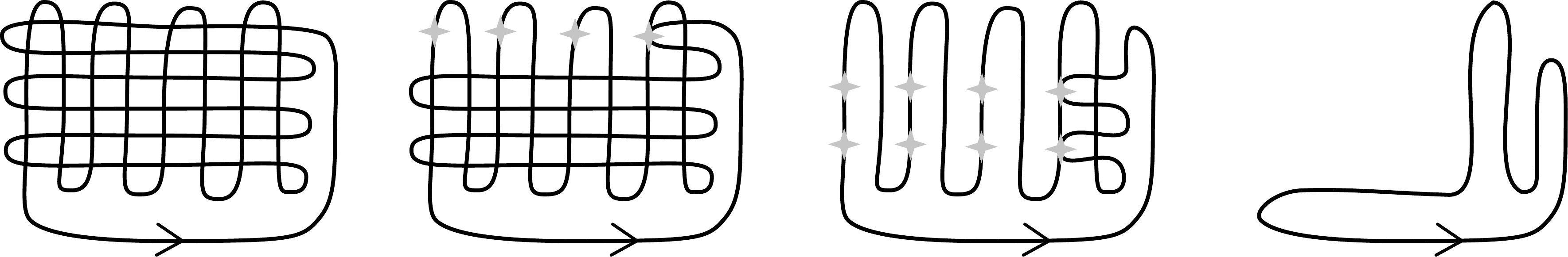}
	\caption{}
	\label{fig:sol11u}
\end{figure}

Immersion $K^u$ is regularly homotopic to the standard curve $K_1 ,$ a circle, through $12$ negative dst, see Figure~\ref{fig:sol11u}, so
$$J^+(K^u) = J^+(K_1) + 2 \cdot 12 = 24 .$$

Immersion $K^v$ is regularly homotopic to the standard curve $K_0 ,$ a figure eight, through $12$ negative dst, similar to $K^u ,$ so
$$J^+(K^v) = J^+(K_0) + 2 \cdot 12 = 24 .$$

\subsubsection*{12. Proving Viro's Formula (page \pageref{ex:12})}

Viro's formula can be proven in \emph{exactly} the same way as the proof we did for Proposition~\ref{prop:rotequation} [Rotation number from winding numbers]. The full proof in this fashion can be found on page 7 of \emph{The $J^{2+}$-Invariant for Pairs of Generic Immersions} by \emph{Hanna Häußler} \cite{hanna:paper}.

\newpage

\pagestyle{litpage}
\nocite{*}
\section*{References}
\addcontentsline{toc}{section}{\hspace{1.4em}References}
\printbibliography[
heading=none,
title={References}
]

\end{document}